\documentclass[11pt]{amsart}
\usepackage{amssymb}
\def\cal{\mathcal}
\def\Bbb{\mathbb}
\def\frak{\mathfrak}

\swapnumbers

\newtheorem{theorem}[subsection]{Theorem}
\newtheorem{lemma}[subsection]{Lemma}

\newtheorem{definition}[subsection]{Definition}

\newtheorem{remark}[subsection]{Remark}

\newcommand{\btimes}                          
{\mathop{\raise1pt\hbox{$\raise-1pt\hbox
to0pt{\small\hskip0.5pt$\square$\hss}\times$}}}
\newcommand{\sfrac}[2]{{\textstyle \frac{#1}{#2}}}
\newcommand{\tcr}{T^{\operatorname{CR}}}

\newcommand{\bbc}{{\Bbb C}}
\newcommand{\bbp}{{\Bbb P}}

\newcommand{\cg}{{\cal G}}

\newcommand{\cl}{{\cal L}}

\newcommand{\ff}[1]{{\frak #1}}
\newcommand{\fg}{{\frak g}}
\newcommand{\fgl}{{\frak g}^L}
\newcommand{\fgr}{{\frak g}^R}

\newcommand{\fsl}{\ff s\ff l}
\newcommand{\fsu}{\ff s\ff u}

\newcommand{\cx}{{\cal X}}

\newcommand{\SL}{\operatorname{SL}}
\newcommand{\SU}{\operatorname{SU}}

\newcommand{\om}{\omega}
\newcommand{\ze}{\zeta}
\newcommand{\et}{\eta}
\newcommand{\ka}{\kappa}
\newcommand{\si}{\sigma}

\newcommand{\La}{\Lambda}
\newcommand{\ph}{\varphi}

\newcommand{\al}{\alpha}
\newcommand{\be}{\beta}
\def\o{\circ}
\newcommand{\x}{\times}
\newcommand{\omi}{\om^{-1}}
\newcommand{\C}{{\mathbb C}}
\renewcommand{\Re}{\operatorname{Re}}

\newcommand{\lra}{\leftrightarrow}

\newcommand{\dyn}[2]{
\begin{picture}(28,15)(0,0)
\put(5,0){\hbox to0pt{\hss$\times$\hss}}
\put(5,8){\hbox to0pt{\hss\tiny $#1$\hss}}
\put(5,2.7){\line(1,0){18}}
\put(23,0){\hbox to0pt{\hss$\times$\hss}}
\put(23,8){\hbox to 0pt{\hss\tiny $#2$\hss}}
\end{picture}}

\def\sideremark#1{\ifvmode\leavevmode\fi\vadjust{
\vbox to0pt{\hbox to 0pt{\hskip\hsize\hskip1em
\vbox{\hsize3cm\tiny\raggedright\pretolerance10000
\noindent #1\hfill}\hss}\vbox to8pt{\vfil}\vss}}}

\begin{document}
\headheight=6mm
\hoffset=-5mm
\title[Hyperbolic and Elliptic CR-manifolds of codimension two]{The Geometry
of Hyperbolic and Elliptic CR-manifolds of codimension two}

\author{Gerd Schmalz and Jan Slov\'ak}

\date{}

\begin{abstract}
The general theory of parabolic geometries is applied to the study of
the normal Cartan connections for all hyperbolic and elliptic
6-dimensional CR-manifolds of codimension two.  The geometric meaning of the
individual components of the torsion is explained and the chains of
dimensions one and two are discussed. 
\end{abstract}
  
\maketitle

\frenchspacing

There have been many attempts to use some ideas going back up to Cartan, in
order to understand the geometry of CR-manifolds. In the codimension one
cases, the satisfactory solution had been worked out in the seventies, see
\cite{Ta0, Ta1, CheM}, but the higher codimensions have not been understood
yet in a comparable extent. In this paper, the recent general theory of the
so called parabolic geometries is applied. In particular, we use the
approach developed in \cite{CSch, Sl2}, see also \cite{Ta, Yam} for earlier
results. Relying on recent achievements by the authors, a clean and quite
simple construction of the normal Cartan connection is presented. This
Cartan connection replaces the absolute parallelisms from \cite{EIS} by more
powerful geometric tools and it enables the detailed study of geometrical
and analytical properties of the CR structures. Consequently the resulting
geometric picture is much more transparent and surprising new results are
obtained.

The main advantage of our approach is the fully coordinate-free handling of
the normal Cartan connection and its curvature. Thus we are able to
translate the cohomological properties of the structure algebras into full
geometrical understanding of the curvature obstruction, without writing down
the curvature components explicitly. The initial section introduces the CR
structures and provides a brief exposition of distinguished second order
osculations of the surfaces by quadrics. Then we observe, that this
osculation transfers enough data from the quadric to apply the general
construction of normal Cartan connections, due to \cite{Ta, CSch}. This
leads easily to the main Theorems \ref{main-hyperbolic} and
\ref{main-elliptic}. In fact, the Cartan connections are constructed also
for certain abstract CR-manifolds though the embedded ones have many
distinguished properties.  The third section is devoted to the exposition of
the generalities on parabolic geometries modelled over $|2|$-graded algebras
and provides the proof of the existence of the normal connections.

Next we study the local geometry of the
hyperbolic points in detail. We recover easily all known facts from
\cite{EIS}, but we go much further. In particular, we identify the complete 
geometric obstructions against the integrability of the almost product 
structure on the tangent bundle (Theorem \ref{torsion-product1}), the
integrability of the almost complex structure on the tangent CR space
(Theorem \ref{torsion-J1}), and
the compatibility of the almost product and almost complex structures
(\ref{torsion-rest}). 
It turns out that the latter two
obstructions always vanish on the embedded hyperbolic CR-structures which
results in automatic vanishing of several algebraic brackets. In particular,
the whole hyperbolic CR-manifold $M\subset {\Bbb C}^4$ is a product of two
3-dimensional CR-manifolds if an only if its almost product structure is
integrable, see Theorem \ref{hyperbolic-torsion-free}. 
Finally we discuss the chains of dimensions one and two. 

Following our intuition, the geometric properties at hyperbolic points
have been expected to have their counterparts in the local geometry at the
elliptic points, cf. remarks and open problems in \cite{EIS}.  This is the
subject of Section \ref{elliptic}. In particular, we observe that the roles
of almost complex and almost product structures are swapped. Thus, there is
an almost complex structure on the whole tangent bundle $TM$ and we
distinguish the algebraic brackets obstructing its integrability in Theorem
\ref{torsion-J2}.  The obstructions against the integrability of the almost
product structures on the complex spaces $T^{CR}M$ and their compatibility
with the almost complex structures vanish automatically for the embedded
elliptic CR-structures. They are discussed in \ref{torsion-product2},
\ref{torsion-rest2}. The analogy to the product property of torsion-free
hyperbolic geometries is the holomorphic normal Cartan connection in the
elliptic case, see Theorem
\ref{torsion-free2}. Finally we prove that for torsion-free elliptic
geometries, there are unique one-dimensional complex chains in all
complex directions transversal to the complex
subbundle $T^{CR}M$ (Theorem \ref{chains-elliptic}).

The last section collects some conclusions and remarks on future
applications. The necessary cohomologies are computed in Appendix \ref{B}
while some more details on
the normalized osculations and the discussion of chains on the hyperbolic
and elliptic quadrics is postponed to Appendix \ref{C}. 

The whole paper stresses the differential-geometric properties and we have
confined the analytical problems and consequences to a few remarks. The
function theoretical aspects will be discussed elsewhere. 

\smallskip
\noindent{\bf Acknowledgements.} The authors should like to mention very
helpful discussions with Andreas \v Cap and Vladim\'{\i}r Sou\v cek. The
research originates in discussions during the stay of both authors at
University of Adelaide. Further support by SFB 256 at Universit\"at Bonn, 
Masaryk University in Brno, and GA\v CR Grant Nr. 201/96/0079 has 
been important. 
Essential part of the writing was undertaken during the stay of the second
author at the Max Planck Institute for Mathematics in Sciences in Leipzig.

\section{CR-structures of codimension two}\label{intro}

Let $M$ be a real submanifold in the complex space $\bbc^N$. Then there is
the CR-subbundle $\tcr M= TM\cap J(TM)$ which consists of all vectors
$\xi_x\in T_xM$ such that the canonical complex structure $J$ on $\bbc ^N$
maps $\xi_x$ to  $J(\xi_x)\in T_xM$. We say that the 
{\em CR-codimension of $M$}
is $k$ if $\operatorname{dim}M$ is $2n+k$ and $\operatorname{dim}\tcr M$ is
$2n$. By means of the implicit function theorem, we may use a
holomorphic projection of $\bbc^N\to \bbc^{n+k}$ and express $M$ locally as
$$
\operatorname{Im} w_\nu = f_\nu(z,\bar z, \operatorname{Re}w),\quad
\nu=1,\dots,k
$$
where $z=(z_1,\dots,z_n)$, $w=(w_1=u_1+iv_1,\dots,w_k=u_k+iv_k)$ are
coordinates in $\bbc^{n+k}$ and $f(0)=0$, $df(0)=0$. Geometrically this means
that the origin belongs to $M$ and $T_0 M$ is just $\{v=0\}$. By means of 
further biholomorphic transformation of second order we are able to 
eliminate the ``harmonic'' part of the second order term in $f$:
$$\Re \sum \frac{\partial^2 f}{\partial z_i\partial z_j}_{|0} z_i z_j + 
2\Re \sum \frac{\partial^2 f}{\partial z_i\partial u_j}_{|0} z_i u_j 
+ \frac{1}{2} \sum \frac{\partial^2 f}{\partial u_i\partial u_j}_{|0} 
u_i u_j.$$    
Only the hermitian part in the second order term of $f$:
$$h(z,\bar{z})=\frac{1}{2}\sum \frac{\partial^2 f}{\partial z_i\partial 
\bar{z}_j}_{|0} z_i \bar{z}_j$$ 
will remain, thus we achieve that $M$ is given by
\begin{equation}\label{osculation}
v=h(z,\bar z) + \operatorname{O}(3)
\end{equation} 
at a neighborhood of the origin. For more details see \cite{Sch}.  
The vector-valued hermitian form
$h$  shall be denoted by $\langle z,z\rangle$ in the sequel. 
The submanifold $M$ is called
Levi non-degenerate (at the origin) if the scalar components of $\langle
z,z\rangle$ are linearly independent and do not have a common annihilator.
The Levi form $\frac{1}{2i}h$ is given by means of the standard Lie bracket 
$\{\ ,\ \}_{\text{Lie}}$ of vector fields modulo the
complex subspace $T^{CR}M$, 
$\xi\mapsto\{\xi_x,J\xi_x\}_{\text{Lie}}\in TM/\tcr M$
for the CR-vector fields $\xi,J\xi:M\to \tcr M$. The latter bracket is
algebraic since the standard Lie bracket
composed with the projection onto the quotient is clearly linear over
functions. 

The geometric meaning of (\ref{osculation}) is that $M$  osculates the quadric
$$Q:\,v=\langle z,z\rangle$$
in second order. Both $M$ and $Q$ share the same tangent space,  
CR-tangent space and Levi form at the origin.

Now, let us assume that $M\subset\bbc^4$ is of CR-codimension 2 and assume
further that $M$ is Levi non-degenerate. Thus $M$ is
a smooth real 6-dimensional manifold.

The quadric $Q$  can be always understood as an open domain in the
homogeneous space $G/P$ where $G$ is the group of the automorphisms of the
hermitian form and $P$ its isotropic subgroup of the origin. This means that
the tangent space in the origin carries the $P$-module structure of $\fg/\ff
p$ in a canonical way and the second order data that are $P$-invariant can
be carried over to $M$ from $\fg/\ff p$ to the individual tangent spaces of
$M$ by means of the osculation.
 
Thus, in order to try to study the geometry of $M$ in the spirit of the
general theory as briefly reviewed in Section \ref{A}, we have to
distinguish the possible non-degenerate $\bbc^2$-valued hermitian forms by a
suitable normalization and to analyze the remaining freedom in the
osculation. This has been done in \cite{Lob, Sch}, see Appendix
\ref{C} for a review. In particular, we can achieve one of the following 
three forms for $h(z,\bar z)=\langle z,z\rangle\in
\bbc^2$ by a linear transformation in $z$'s and $v$'s 
\begin{alignat}{2}\label{hyperbolic-quadric}
h^1(z,\bar z)&= z_1\bar z_1,&\quad h^2(z,\bar z)&= z_2\bar z_2\\
\label{parabolic-quadric}
h^1(z,\bar z)&= z_1\bar z_1,&\quad h^2(z,\bar z)&= \operatorname{Re}z_1\bar z_2\\
\label{elliptic-quadric}
h^1(z,\bar z)&= \operatorname{Re}z_1\bar z_2,&\quad h^2(z,\bar z)&=
\operatorname{Im}z_1\bar z_2
\end{alignat}
and we refer to these cases as to hyperbolic, parabolic, and elliptic,
respectively. The normalization (\ref{osculation}) with one of these 
hermitian forms $h$ is given uniquely up to the isotropic subgroup of the
origin in the group of all biholomorphic automorphisms of $Q\subset{\Bbb
C}^4$.

We say that a point $x\in M$ is {\em hyperbolic} or {\em parabolic} or {\em 
elliptic} if the
osculating quadric at $x$ is of that type. Apparently, the set of all
hyperbolic points is open and the same for the elliptic ones. The
CR-structure on $M$ is called {\em hyperbolic}, or {\em parabolic}, or 
{\em elliptic}, if all points of $M$ are of the same type.

Let $M\subset \bbc^4$ be a CR-structure of codimension two, such that all
its points are either hyperbolic or elliptic. As discussed above, the choice
of the canonical form of the osculating quadric $Q=G/P$ reduces the freedom
in the osculation (\ref{osculation}) to the isotropic subgroup of the origin
in $G/P$ and this allows to transfer the $P$-invariant data of first and
second order from the origin of
$Q$ to the individual tangent spaces in all points of $M$. 

We recall the details on the resulting groups 
\begin{equation}\label{hyperbolic_group}
G=\bigl((\SU(2,1)/{\Bbb Z}_3)\times (\SU(2,1)/{\Bbb Z}_3)\bigr) 
\rtimes {\Bbb Z}_2
\end{equation}
in the hyperbolic case, and
\begin{equation}\label{elliptic_group}
G=(\SL(3,{\Bbb C})/{\Bbb Z}_3)\rtimes {\Bbb Z}_2
\end{equation}
in the elliptic case, $P$, $G_0$, and their Lie algebras 
in Appendix \ref{C}. At the moment, let us notice that in both cases the Lie
algebra $\fg$
carries the $|2|$-grading
$\fg=\fg_{-2}\oplus\fg_{-1}\oplus\fg_0\oplus\fg_1\oplus\fg_2$, $\ff
p=\fg_0\oplus\fg_1\oplus\fg_2$, and the
subgroups $P$ and $G_0$ have all properties discussed in
\ref{graded-algebras} below.
In particular $P$ is the subgroup of all elements whose adjoint action
leaves the $\ff p$-submodules in $\fg$ invariant, while
$G_0$ consists of all elements which leave the components $\fg_i$ invariant.
Thus, the tangent space $T_xM$
at each point $x\in M$  is identified with the $P$-module  $\ff
g/\ff p$ which is the tangent space to the 
osculating quadric $Q$ at its origin, the normalized osculation
transfers the $P$-submodule $\fg_{-1}\subset \fg/\ff p$ to $T^{CR}M\subset
TM$, and the algebraic structure of 
$\fg/\ff p$ is carried over to the associated
graded tangent space $\operatorname{Gr}TM=(TM/T^{CR}M)\oplus T^{CR}M$. 

\begin{lemma}\label{algebraic-brackets}
Let $M\subset {\Bbb C}^4$ be a hyperbolic or elliptic 6-dimensional 
CR-manifold. Then all algebraic brackets 
$T^{CR}M\x T^{CR}M\to T^{CR}M$ and 
$T^{CR}M\x T^{CR}M\to TM/T^{CR}M$ on the real graded tangent space
$\operatorname{Gr}TM$, and the
analogous algebraic brackets on the complexified graded tangent space
$\operatorname{Gr}T_{\Bbb C}M$
are obtained via the osculation from the corresponding brackets at the
origin of the quadric. 

In particular, the algebraic Lie bracket $\{\ ,\
\}_{\text{Lie}}$ on $\operatorname{Gr}TM$ 
coincides with the algebraic bracket carried over by the osculation
(\ref{osculation}). 
\end{lemma}  
\begin{proof}
The Lie bracket on $\fg_-=\fg_{-2}\oplus\fg_{-1}$ 
is $G_0$-equivariant, and so the osculation
(\ref{osculation}) induces an algebraic bracket on the associated graded
vector bundle $\operatorname{Gr}TM$.
A neighborhood of the origin in $Q$ can be identified with the exponential
image of $\fg_-$ in $G$ and the Lie bracket in $\fg_-$ is given by the 
usual Lie brackets of the left invariant vector fields on $G$. By means of
the osculation, we can project these fields onto $M$
locally and clearly the algebraic bracket $T^{CR}M\x T^{CR}M\to TM/T^{CR}M$ 
induced by the Lie bracket of
vector fields on $M$ coincides with that one carried over from $\fg_-$ 
by the osculation. Obviously, the result will not be effected by the action of
an element in $P_+$ on $Q$ (i.e. by the possible change of the osculation).

All other algebraic brackets on the real graded tangent space can be treated
in exactly the same way, provided they are $P$-invariant on the quadric. 
The action of an element of $G_0$ always 
commutes with the osculation while the
action of $P_+$ is not visible in all our cases. Indeed, the action is
trivial if all arguments and values are in $T^{CR}M$, while the
contributions of the action is factored out in the case of the brackets
$T^{CR}M\x T^{CR}M\to TM/T^{CR}M$. Similarly, the left invariant 
vector fields in the complexified tangent spaces on the quadric 
can be mapped into complex vector fields on $M$ and the above arguments
apply as well.

Let us notice, however, that the possible algebraic brackets taking
some arguments in $TQ/T^{CR}Q$ are never $P$-invariant.
\end{proof}

The latter lemma turns out to be the most crucial point for our further
development. Indeed, there is the general theory of the so called parabolic
geometries which we adapt for our purposes in the next section. In
particular, Theorem \ref{A-main} due to \cite{Ta,CSch} will provide the
canonical principal bundles together with canonical Cartan connections for
all hyperbolic and elliptic 6-dimensional CR-manifolds with CR-codimension
two (see the beginning of Section \ref{A} for definitions and more
explanation).  We should also like to mention already now that the complete
proof of Theorem \ref{A-main} is in fact constructive, it is based on well
known facts from representation theory, and it is even shorter and simpler
than the ad hoc construction of the absolute parallelisms in \cite{EIS}. 
The ultimate results read as follows:

\begin{theorem}\label{main-hyperbolic} 
On each 6-dimensional hyperbolic CR-manifold $M\subset \bbc^4$ of 
CR-codimension two, there
is the unique normal Cartan connection $\om$ of type $(G/P)$ 
on the principal fibre bundle $\cg\to M$, up to isomorphisms. The subgroup
$P$ is the subgroup of all elements in $G$ from (\ref{hyperbolic_group}) 
which respect
the $\ff p$-module filtration on $\ff s\ff u(2,1)\oplus \ff s\ff u(2,1)$.
\end{theorem}

\begin{theorem}\label{main-elliptic} 
On each 6-dimensional elliptic CR-manifold $M\subset \bbc^4$ of 
CR-codimension two, there
is the unique normal Cartan connection $\om$ of type $(G/P)$ 
on the principal fibre bundle $\cg\to M$, up to isomorphisms. The subgroup
$P$ is the subgroup of all elements in $G$ from (\ref{elliptic_group}) which respect
the $\ff p$-module filtration on $\ff s\ff l(3,{\Bbb C})$.
\end{theorem}

For the proof of these theorems see \ref{main-proofs} below.
The reason why the methods of \cite{EIS} could
not produce a principal fibre bundle $\cg$ with structure group $P$ and a 
normal Cartan connection on $\cg$, was hidden in the initial choice of the
normalization which had to produce a Cartan connection without torsion. In
our approach, the torsions are the important parts of the curvature which
are easily observable on the CR-manifold itself. The Sections
\ref{hyperbolic} and \ref{elliptic} are
basically dealing with the consequences of the vanishing of the individual
components of the torsion of the canonical Cartan connection for 
the hyperbolic and elliptic local geometries.

\section{Parabolic geometries modelled over $|2|$-graded Lie algebras}
\label{A}

The aim of this section is to introduce the reader to the so called parabolic
geometries, but we shall concentrate on the cases similar to the real forms 
of the two-graded complex Lie algebra $\ff g=\fsl(3,\bbc)\oplus \fsl(3,\bbc)$.
Beside well known facts, we shall also have to adapt and extend some
points.

Let us notice first that the general ideas go back to E.~Cartan and his
notion of ``espace generalis\'e''. The interest in the parabolic structure
groups was pointed out by Fefferman, \cite{Fef}, in connection with problems
in conformal and CR geometries. Extensive study was undertaken even earlier
by Tanaka (see \cite{Ta} and the references therein), motivated by a class
of equivalence problems for differential systems. Tanaka's approach was
developed further, see e.g. \cite{Mor,Yam}. Motivation coming from twistor
calculus led to another direction of related research, see e.g.
\cite{BaE,BEG,Gosrni}.
General background and an introduction to the subject may be also found in
\cite{Sha}. The exposition in this section extends the development in
\cite{CSS1, CSch} and follows mainly \cite{Sl2}.
 
\subsection{Graded Lie algebras}\label{graded-algebras}
Let $\ff g$ be a $|2|$-graded Lie algebra, $\ff p$ and $\ff g_-$ its 
subalgebras:
\begin{gather*}
\ff g=\ff g_{-2}\oplus \fg_{-1}\oplus\ff g_0\oplus\ff
g_1\oplus\ff g_{2}\\
\ff g_- = \ff g_{-2}\oplus\ff g_{-1},\quad
\ff p = \ff g_0\oplus\ff
g_1\oplus\ff g_{2}
.\end{gather*}
Further, let $G$ be a Lie group with the Lie algebra $\ff g$. Then there is
the subgroup $P\subset G$ of elements whose adjoint representations on
$\ff g$ preserve the filtration by $\ff p$-submodules $\ff
g_i\oplus\fg_{i+1}\oplus\dots \oplus\fg_2$ and there also is the subgroup
$G_0\subset P$ of all elements whose adjoint representation leaves invariant
all $\fg_i$. Thus the components $\fg_i$ of the grading can be understood 
as $G_0$-submodules, but also as the factors in the graded $P$-module
components associated to the $P$-module filtration. Similarly we define the
$|k|$-graded algebras $\fg= \fg_{-k}\oplus\dots\oplus\fg_k$. 

In the sequel, we shall deal with semi-simple $|2|$-graded Lie algebras
exclusively. It is well known that all graded semi-simple Lie algebras are
sums of $|k|$-graded algebras for suitable $k$'s and the subgroups $P\subset
G$ are always suitable real forms of parabolic subgroups $P_{\Bbb C}\subset
G_{\Bbb C}$ in the complexification. The exposition below extends
easily to general semi-simple
$|k|$-graded Lie algebras and the corresponding parabolic structures, as
discussed in \cite{CSch} for example. Many geometric and algebraic
properties of these geometries are deduced in \cite{CSS4}.

\subsection{Cartan connections} 
The homogeneous space $p:G\to G/P$
is equipped with the left Maurer-Cartan form $\om\in C^\infty(\cg, \ff g)$.
This is the prototype of a geometry modelled over the homogeneous space
$G/P$.
In general, a {\em Cartan geometry of type $G/P$} is 
a principal fibre bundle $p:\cg\to M$ over a smooth manifold $M$, equipped
with a $\ff g$-valued one-form $\om\in C^\infty(\cg,\ff g)$ satisfying
\begin{itemize}
\item $\om(\ze_X(u))=X$ for all $X\in \ff p$ and fundamental fields
$\ze_X$ on $\cg$,
\item $\om$ is right-invariant, i.e. $(r^b)^*\om =
\operatorname{Ad}(b^{-1})\o\om$ for all $b\in P$,
\item the restrictions $\om_{|T_u\cg}: T_u\cg\to \ff g$ are linear
isomorphisms, i.e. the obvious 
mapping $T\cg\to \cg\x \ff g$ is
a diffeomorphism.
\end{itemize}

The {\em homomorphisms} of Cartan geometries are those principal fibre 
bundle morphisms which respect the Cartan
connections. The {\em flat} Cartan geometry is the homogeneous space $G\to
G/P$ with the Maurer-Cartan form $\om$.

Let us also observe that the above absolute parallelisms $\om$ turn out to be
special cases of principal connections $\tilde \om$ on the principal bundle
$\tilde \cg=\cg\x_P G$ with structure group $G$. Indeed, the connection
forms of all principal connections on $\tilde G$ whose horizontal
distributions do not meet the tangent space $T\cg\subset T\tilde \cg$
restrict to forms $\om$ with the required properties. See e.g. \cite{Sha}
for more comments. 

\subsection{Normal coordinates}
For each $X\in \ff g$, the parallelism $\om$ defines the vector field
$\om^{-1}(X)$ on $\cg$. The {\em horizontal vector fields $\om^{-1}(X)$} 
on $\cg$ are those with $X\in \ff g_-$ and their values span the {\em
horizontal distribution} on $\cg$. Due to the third property of $\om$, the
choice of a frame $u\in\cg$ defines an injective smooth mapping of a
neighborhood of zero in $\fg_-$
\begin{equation}
\fg_-\ni X\mapsto \operatorname{Fl}^{\om^{-1}(X)}_1(u)
\end{equation}
defined by means of the flows of the vector fields $\om^{-1}(X)$.
The tangent space of its image at $u$ belongs to the horizontal distribution
on $\cg$ and its composition with the projection $p:\cg\to M$ defines the
locally defined mapping 
\begin{equation}\label{A-coordinates}
\tilde u:\fg_-\to M,\quad X\mapsto p(\operatorname{Fl}^{\om^{-1}(X)}_1(u))
\end{equation}
which is diffeomorphic on a neighborhood of the origin. We call $\tilde u$
the {\em normal coordinates} on $M$ given by the frame $u$.
At the same time, 
$\tilde
u$ induces the local trivialization $\si_u$, 
\begin{equation}\label{A-trivialization}
\si_u:M\to \cg,\qquad M\ni \tilde u(Y)\mapsto
\operatorname{Fl}^{\om^{-1}(Y)}_1(u)\in \cg.
\end{equation}
Clearly, the normal coordinates around a fixed
point $x\in M$ are parameterized by elements in $P$ and they generalize the
usual normal coordinates of affine connections on manifolds. 
The general concept of the normal coordinates has been introduced 
and studied in \cite{Sl2}.

\subsection{Chains}\label{chains-general}
The notion of normal coordinates suggests a straightforward generalization of
the geodetical curves. For each choice of $X\in \fg_{-}$ and $u\in\cg$ we
define the {\em 1-chain} $\al^{u,X}:{\Bbb R}\to M$ on a neighborhood of
$0\in{\Bbb R}$ by
$$
\al^{u,X}(t)=p(\operatorname{Fl}^{\omi(X)}_t(u))
.$$
Clearly the tangent direction to the 1-chain $\al^{u,X}$ at its origin is
the vector $Tp.\omi(X)(u)$ since the tangent bundle $TM$ is identified with 
the associated bundle $\cg\x_P (\fg/\ff p)$ via the adjoint representation,
$\{u,X\}\mapsto Tp.\omi(X)$. 
In particular we see immediately that many different
1-chains may share the same tangent direction.

The 1-chains have been studied under various names like Cartan's
circles or generalized circles, see e.g. \cite{Sha}, and the chains
introduced by Chern and Moser for CR-geometries of codimension one are
exactly the 1-chains with $X\in \ff g_{-2}$. Since
$\operatorname{dim}\fg_{-2}=1$ for these geometries, the latter 1-chains 
coincide with the chains defined below. 

All 1-chains corresponding to a fixed frame $u$ yield exactly the normal
coordinates with origin at $p(u)$ and the transformation rules for these
coordinates under the change of $u$ may be quite complicated, in general. On
the other hand, the 1-chains corresponding to the parameters $\{u,X\}$
with $X\in\ff g_{-2}$ have very specific properties. We define the {\em
chain} $\be^u:\ff g_{-2}\to M$ by the formula
$$
\be^u(X)=p(\operatorname{Fl}^{\omi(X)}_1(u))
.$$
Thus the chains are parameterized submanifolds in $M$ of dimension
$\operatorname{dim}\fg_{-2}$.

\begin{remark}\label{chains-structure}\em
The importance of the chains grows whenever they are given uniquely by
their tangent directions in the origin. Another important question is
whether two different chains may intersect nontrivially in each small
neighborhood of their common origin. The answer to these questions is usually
very easy because of the following equivalent definition of 1-chains by
means of their developments into the associated bundle
$FM=\tilde\cg\x_G(G/P)$. 

The principal connection $\tilde \om$ on $\tilde \cg$ provides the induced
(generalized) connection on the bundle $FM$ and there is the canonical
embedding of $M$ into $FM$, $p(u)\mapsto \{u,[e]\}$. Thus each curve
$\al(t)\in M$ can be mapped by the parallel
transport of $\tilde\om$ into a curve $\tilde\al$ in the fibre over 
$\al(0)$. This curve $\tilde \al$ is called the {\em development of the
curve $\al$}.  Clearly, the germs of curves through $\al(0)$ are in bijective
correspondence with the germs of  their developments.

Now, our definition of the 1-chains can be easily rephrased as follows. The
1-chains are exactly the curves $\al$ whose developments $\tilde \al$ are
given by one-parametric subgroups in $G$, i.e.
$\tilde\al^{u,X}=\{u,[\operatorname{exp}tX]\}$. See e.g. \cite{Sl2} for
more details. 

Since our chains $\be^u$ are obtained via 1-chains, all
structural questions mentioned above are obtained by the discussion of the
chains in the homogeneous case.   
\end{remark}

\subsection{Curvature and torsion}
The structure equation
$$
d\om= -\frac12[\om,\om]+K
$$
defines the $\ff g$-valued horizontal 2-form $K\in\Omega^2(\cg,\ff g)$. If we evaluate the structure equation on
two horizontal vector fields we obtain the so called frame form of the
curvature, the equivariant 
function $\ka\in C^{\infty}(\cg, \Lambda^2\ff g_-^*\otimes\ff g)^P$
\begin{align*}
\ka(u)(X,Y)&=K(\om^{-1}(X),\om^{-1}(Y))(u)
\\
&=[X,Y]-
\om([\om^{-1}(X),\om^{-1}(Y)](u))
.\end{align*}
The Cartan geometry is locally isomorphic to the flat one if and only if its
curvature vanishes. 

If $\ff g$ is semi-simple, then $P$ is a parabolic subgroup of the
semi-simple group $G$ and we then refer 
to the above geometries as to {\em parabolic geometries 
of type $G/P$}. 

The curvature $\ka$ has values in the space of cochains of the Lie algebra
cohomology $H^*(\ff g_-,\ff g)$. The grading on $\ff g$ induces the grading on the space of cochains. The
homogeneous cochains of degree $k$ are those which map $\ff g_i\wedge \ff
g_j$ into $\ff g_{i+j+k}$ and this grading is respected by the Lie algebra
cohomology differential  $\partial$. For each cochain $\et\in
\La^k\fg_{-}\otimes A$ with values in a $\fg_-$-module $A$ the differential
is given by
\begin{equation}\label{partial}
\begin{aligned}
\partial \et(X_0,\ldots,X_k) =\ &\sum_{i=0}^k (-1)^i
X_i.\et(X_0,\hat{\ldots},X_k) \\&+
\sum_{i<j}(-1)^{i+j}\et([X_i,X_j],X_0,\hat{\ldots},X_k)
\end{aligned}
\end{equation}
where the dot in the first summand means the $\fg_{-}$-module action while
the hats denote the obvious omitions.

In particular, the whole curvature splits into the
homogeneous parts $\ka^{(k)}$
$$
\ka=\sum_{k=-\ell+2}^{3\ell}\ka^{(k)}
$$
where $\ell=2$ is the length of the grading. 
On the other hand, we may split $\ka$ according to its values. In
particular, there is the {\em torsion part} $\ka_-$ with values in $\ff g_-$
$$
\ka=\sum_{i=-\ell}^{\ell}\ka_{i}\qquad \ka_-=\ka_{-\ell}\oplus\dots\oplus
\ka_{-1}\qquad\ka_{\ff p}=\ka_0\oplus\dots\oplus\ka_{\ell}
.$$
The torsion has a simple geometrical meaning: Let us define the horizontal
bracket $[\xi,\et]_h$ on the space $\cx_h(\cg)$ 
of all fields belonging to the horizontal
distribution on $\cg$ by the standard Lie bracket followed by horizontal
projection. By the very definition, the torsion of $\om$ vanishes if and
only if  the mapping $\ff g_-\ni X\mapsto \om^{-1}(X)\in \cx_h(\cg)$ is a
Lie algebra homomorphism. 

\subsection{Regular and normal connections}\label{Hodge-structure}
We say that the parabolic geometry $(\cg,\om)$ 
is {\em regular} if $\ka^{(k)}=0$ for all
$k\le 0$.

In the sequel, we shall always assume $\ff g$ is
semi-simple. Then there is the adjoint of the Lie
algebra cohomology  differential
$\partial$, the codifferential $\partial ^*:
\Lambda^k\ff g_-^*\otimes \ff g\to \Lambda^{k-1}\ff g_-^*\otimes 
\ff g$
.

We say that $\om$ is a {\em normal Cartan
connection} if its curvature is co-closed, i.e. 
$$
\partial^*\o\ka=0\in C^\infty(\cg,\ff g_-^*\otimes \ff g).
$$
Let us recall, that the whole space of cochains decomposes into a sum of
irreducible components as a $\ff g_0$-module. Each such component is either
in the image of $\partial$ or in the image of $\partial^*$ or in the kernel
of both. The latter components are called harmonic and they are in bijective
correspondence with the non-zero cohomologies $H^*(\ff g_-,\ff g)$. 

\begin{theorem}\label{A-curvature} {\em (\cite{Ta,Yam,CSch})}
Let $(\cg,\om)$ be a normal Cartan connection and assume that all components
$\ka^{(j)}$, $j<k$, vanish. Then $\partial\o \ka^{(k)}$ vanishes and so 
all non-trivial irreducible components of
$\ka^{(k)}$ are harmonic.

In particular, the whole curvature of $\om$ vanishes if and only if its
harmonic part does.
\end{theorem}

The latter theorem is a straightforward consequence of the important 
Bianchi identity for Cartan geometries:
\begin{equation}\label{A-Bianchi}
\begin{aligned}
\partial \ka^{(k)}(X,Y,Z) = &-\sum_{\operatorname{cyclic}}
\sum_{i=1}^{k-1}\ka^{(k-i)}
(\ka^{(i)}(X,Y),Z) 
\\ &-\sum_{\operatorname{cyclic}}\cl_{\om^{-1}(Z)}\ka^{(k+|Z|)}(X,Y)\bigr)
\end{aligned}
\end{equation}
where the sum is the cyclic sum over $X,Y,Z\in \ff g_-$, and $|Z|=j$ if
$Z\in\ff g_j$. See e.g. \cite{CSch} for more details.

\subsection{The underlying geometry on $M$}
A part of the Cartan geometry $(\cg, \om)$ is visible directly on the
underlying manifold $M$ and, fortunately, these
data are sufficient in order to reconstruct the Cartan connection
completely. This is the core of our approach to the CR structures in this
paper. As before we shall restrict ourselves to the
$|2|$-graded cases below, but the discussion extends easily to the general
case. 

The $P$-module structure on $\ff g$ (defined via the $\operatorname{Ad}$
representation) determines the filtration by $P$-submodules 
\begin{gather*}
\ff g=V^{-2}\supset V^{-1}\supset V^0\supset V^1\supset V^2=\ff g_2\\
V^k=\ff g_k\oplus\dots\oplus\ff g_2\subset \ff g, \quad k=-2,-1,0,1,2
.\end{gather*}
This in turn defines the filtration on $T\cg$ 
\begin{gather*}
T\cg=T^{-2}\cg\supset T^{-1}\cg\supset T^{0}\cg\supset T^{1}\cg\supset
T^2\cg
\\
T^k_u(\cg)=\om^{-1}(u)(V^k),\quad k=-2,-1,0,1,2,\ u\in\cg
.\end{gather*}  
The right invariance of $\om$ yields 
\begin{equation}
\om^{-1}(u.b)(X)=Tr^b.\om^{-1}(u)(\operatorname{Ad}(b).X)
\end{equation}
and so the latter filtration on $\cg$ is $P$-invariant. The $P$-invariant
projection $p:\cg\to M$ defines then the filtration
$$
TM=T^{-2}M\supset T^{-1}M\supset \{0\}
.$$
Moreover, each fixed frame $u\in\cg$ with $p(u)=x\in M$ 
determines the linear isomorphism of
filtered vector spaces 
$$
\hat u:\ff g_-\to T_xM \qquad X\mapsto Tp.\om^{-1}(X)(u)
$$ 
and on the level of the associated graded spaces we obtain 
the linear isomorphism
$$
\hat u:V^{-2}/V^{-1}\oplus V^{-1}/V^{0}\simeq
\ff g_{-2}\oplus\ff g_{-1}\to T_x^{-2}M/T_x^{-1}M\oplus T_x^{-1}M
.$$
The whole structure group $P$ is a semidirect product of its reductive 
subgroup $G_0$ (corresponding to the Lie algebra $\ff g_0$) and the subgroup
$P_+$ which corresponds to $\ff p_+=\ff g_1\oplus\ff g_2$.
Obviously, the latter identification $\hat u$ 
does not change if we replace the frame $u$ by $u.b$ with
$b\in P_+$. Thus we have identified the {\em graded tangent bundle}
$\operatorname{Gr}TM$ with the associated vector bundle to the principal
bundle $\cg_0=\cg/P_+$ whose standard fibre is the $G_0$-module $\ff g_-$. 
In particular the Lie bracket on $\ff g_-$ is transfered to the algebraic
bracket $\{\ ,\ \}_0$ by
$$
\{\xi_x,\et_x\}_0= \hat u([\hat u^{-1}(\xi_x),\hat u^{-1}(\et_x)]),
\quad \xi_x,\et_x\in \operatorname{Gr}T_xM, u\in \cg
.$$
Notice that this definition does not depend on the choice of $u$ since
$\hat u$ is independent of the 
action of $P_+$ and the Lie bracket on $\ff g_-$
is $G_0$-equivariant. Since our $G_0$-structure on $\operatorname{Gr}TM$ is
defined by the Cartan connection, we may choose representing vectors 
$\bar\xi_x\in T^i_xM$,
$\bar\et_x\in T_x^jM$, their covering vectors $\hat \xi_u$,$\hat\et_u\in
T\cg$ and we obtain 
$$
\{\xi_x,\et_x\}_0=\pi(\om^{-1}([\om(\hat \xi_u),\om(\hat \et_u)])(u))
$$
where $\pi$ is the obvious projection $T^{i+j}\cg\to T^{i+j}M\to
T^{i+j}M/T^{i+j+1}M$. 

We shall see in a while that the regular Cartan geometries are
exactly those for which the latter bracket is induced from the Lie bracket
of vector fields in an algebraic way. Since we shall need a good control
over the relations between the brackets of the horizontal vector fields and
some brackets on the underlying manifold in its proof, we shall first
formulate a general lemma based on our concept of the normal
coordinates.

\begin{lemma}\label{A-brackets-lemma}
Let $u\in \cg$ and let $\si_u$ be the corresponding distinguished local 
trivialization of $\cg$, see (\ref{A-trivialization}). Further let
$X,Y\in \fg_-$, and consider the projectable vector fields $\tilde \xi$,
$\tilde\et$ over $M$, such that their restrictions to the image of $\si_u$
coincide with the horizontal fields $\om^{-1}(X)$, $\om^{-1}(Y)$,
respectively. Then $[\om^{-1}(X),\om^{-1}(Y)](u)= 
[\tilde \xi,\tilde\et](u)$. 

Thus, in particular
$$
Tp.[\om^{-1}(X),\om^{-1}(Y)](u)= [Tp.\tilde \xi,Tp.\tilde\et](p(u))
$$
\end{lemma}
\begin{proof} Let us write $\tilde\xi=\om^{-1}(X)+\mu$,
$\tilde\et=\om^{-1}(Y)+\nu$ and compute
their bracket. By the very definition, we obtain 
\begin{align*}
[\tilde\xi,\tilde\et](u)&= \frac d{dt}_{|_0} T(\operatorname{Fl}^{\tilde
\xi}_{-t})\o(\om^{-1}(Y)+\nu)\o (\operatorname{Fl}^{\tilde \xi}_t)(u)
\\
&=\frac d{dt}_{|_0} T(\operatorname{Fl}^{\tilde
\xi}_{-t})\o(\om^{-1}(Y)+\nu)\o (\operatorname{Fl}^{\om^{-1}(X)}_t)(u)
\\ 
&=\frac d{dt}_{|_0} T(\operatorname{Fl}^{\tilde
\xi}_{-t})\o(\om^{-1}(Y))\o (\operatorname{Fl}^{\om^{-1}(X)}_t)(u)
\\
&=[\tilde \xi, \om^{-1}(Y)](u)
\end{align*}
where the first equality follows from the fact that the flows of $\tilde \xi$
and $\om^{-1}(X)$ through $u$ coincide, the next one results from the
vanishing of $\nu$ on the image of $\si_u$. Now, repeating the same
arguments for $[\om^{-1}(Y),\tilde\xi]$, we achieve just the required
equality.
\end{proof}

\begin{lemma}\label{A-brackets}
Let $\om\in\Omega^1(\cg,\ff g)$ be a Cartan connection with a $|2|$-graded
Lie algebra $\fg$. Then
$\ka^{(i)}=0$ for all $i<0$ and the Lie bracket of
vector fields defines an algebraic
bracket $\{\ ,\ \}_{\operatorname{Lie}}$ on the graded vector bundle
$\operatorname{Gr}TM$. 
Moreover, $\ka^{(0)}$ vanishes if and only if the latter
bracket coincides with the algebraic bracket $\{\ ,\ \}_0$ on
$\operatorname{Gr}TM$.
\end{lemma}
\begin{proof}
Recall that the defining equation for
the homogeneous components $\ka^{(k)}(u)(X,Y)$, $k\ne0$, $u\in \cg$, 
$X\in \ff g_i$, $Y\in\ff g_j$ is
\begin{equation}\label{curvature1}
\ka^{(k)}(u)(X,Y)=
-\om_{i+j+k}([\om^{-1}(X),\om^{-1}(Y)](u))
\end{equation}
while the component of degree zero is
\begin{equation}\label{curvature2}
\ka^{(0)}(u)(X,Y)=
[X,Y]-\om_{i+j}([\om^{-1}(X),\om^{-1}(Y)](u))
\end{equation}
Now, consider vector fields $\xi$ in $T^iM$, $\et$ in $T^jM$ and let us
choose elements $X_r\in V^i$, $Y_s\in V^j$ such that
$\xi=Tp.\sum_rf^r\om^{-1}(X_r)$, $\et=Tp.\sum_s g^s\om^{-1}(Y_s)$ with suitable
functions $f^r$, $g^s$ on $\cg$. Then
$$
[\xi,\et] =(Tp. \sum_{r,s}f^rg^s[\om^{-1}(X_r),\om^{-1}(Y_s)])\
\operatorname{mod}T^{i+j+1}M.
$$

The negative homogeneous components $\ka^{(k)}$, $k<0$ have to vanish
because the algebra is $|2|$-graded and so we have no choice of arguments 
for cochains with such homogeneity. The lowest possible case will be a cochain
$\fg_{-1}\x\fg_{-1}\to\fg_{-2}$ of homogeneity zero. The fact that the
Lie bracket of vector fields produces an algebraic bracket on the associated
graded tangent bundle is obvious.  

Now, the two brackets in question may be expressed for all vectors
$\xi_x=\pi(\om^{-1}(X)(u))$ and $\et_x=\pi(\om^{-1}(Y)(u))$ as
\begin{align*}
\{\xi_x,\et_x\}_0&=\pi(\om^{-1}([X,Y])(u))
\\
\{\xi_x,\et_x\}_{\text{Lie}}
&=[Tp.\tilde\xi,Tp.\tilde\et](x)\text{ mod } T^{i+j+1}M
\\&=\pi([\om^{-1}(X),\om^{-1}(Y)](u)\text{ mod }T^{i+j+1}\cg)
\\
&=\pi(\om^{-1}(\om_{i+j}([\om^{-1}(X),\om^{-1}(Y)](u))))
\end{align*}
where $\tilde\xi$ or $\tilde \et$ are some projectable fields 
from the previous Lemma \ref{A-brackets-lemma}.
Thus, according to (\ref{curvature2}), the two brackets equal each other 
if and only if $\ka^{(0)}$ vanishes.
\end{proof}

Now we have got the motivation for the following definition of geometric
structures on manifolds. Let us also remark that the version of the
latter lemma which is valid for all $|k|$-graded structures needs one
more condition. Namely, the existence of the algebraic bracket induced by
the Lie brackets of vector fields is equivalent to the vanishing of all
negative components $\ka^{(k)}$, $k<0$.

\begin{definition}\label{gp-structure} Let $\fg$, $G$, $P$, and $G_0$ be as
in \ref{graded-algebras}. 
A regular $(\ff g, P)$-structure on
a smooth manifold $M$ is a filtration of the tangent bundle $TM$ 
$$
TM=T^{-2}M\supset T^{-1}M
$$
together with the reduction of the structure group of the associated graded
tangent vector bundle $\operatorname{Gr}TM$ to the subgroup $G_0$, such that
the algebraic bracket on $\operatorname{Gr}TM$
induced by the Lie bracket of vector fields coincides with the algebraic Lie
bracket defined by the $G_0$-structure. 
\end{definition}

We may understand the above condition as the requirement that
the subbundle $T^{-1}M$ be reasonably non-involutive. Due to our
restriction to $|2|$-graded algebras we do not need to consider the other 
condition from Lemma \ref{A-brackets} on the Lie brackets
of vector fields, namely that they must not be ``too much non-involutive''.

Surprisingly enough there is the theorem claiming that, apart from a few
exceptions, {\em all\/} regular normal parabolic geometries are uniquely
given by the underlying $(\fg,P)$-structures on the manifolds $M$:

\begin{theorem}\label{A-main}
Let $M$ be a smooth manifold, $\ff g$ a graded semi-simple Lie algebra, $G$
a Lie group with Lie algebra $\ff g$,
and assume that all homogeneous components of the cohomologies
$H^1_\ell(\ff g_-,\ff g)$ with positive degrees $\ell>0$ are trivial. Then
there is a bijective equivalence between isomorphism classes of the regular
$(\ff g,P)$-structures on $M$ and the isomorphism classes of regular
normal Cartan geometries $(\cg,\om)$ over $M$.
\end{theorem}   

For the proof see Section 3 of \cite{CSch}. The computations in
\cite{Yam,CSch}
show that, apart from situations with simple
components in $\ff g_0$, the only exceptions are $\ff g=\fsl(2,\Bbb C)$, 
specific maximal parabolic subalgebras in special linear
algebras in higher dimension ($|1|$-graded examples) and specific maximal
subalgebras in symplectic algebras ($|2|$-graded examples). An equivalent
theorem for the cases $\ff g$ simple and $G$ connected was proved in
\cite{Ta}.

\subsection{Proof of Theorems \ref{main-hyperbolic} and \ref{main-elliptic}}
\label{main-proofs}
The relevant cohomologies for the real forms of $\fsl(3,{\Bbb
C})\oplus\fsl(3,{\Bbb C})$ are computed in Appendix \ref{B}. In particular,
there is no obstruction in the construction of the normal Cartan connections
out of regular $(\fg,P)$-structures according to Theorem \ref{A-main}.
The definition of the relevant $(\fg,P)$-structures by means of the
fundamental second order osculation (\ref{osculation}) was discussed at the
end of Section \ref{intro}, see Lemma \ref{algebraic-brackets}.

\section{The hyperbolic structures}\label{hyperbolic}

In this section, we shall study the consequences of the
algebraic structure of $\ff s\ff u(2,1)\oplus \ff s\ff u(2,1)$ for the
hyperbolic points on 6-dimensional CR-manifolds of CR-codimension 2
$M\subset {\Bbb C}^4$. Thus the Lie groups $G$, $P$, $G_0$, 
as well as the corresponding Lie algebras will be fixed throughout this 
section.

\subsection{Almost product and almost complex structures} As we noticed
already in the proof of Theorem \ref{main-hyperbolic}, there is the relevant 
$(\fg,P)$-structure on $M$.
Since the individual left and right components of
$\fg$ are $P$-submodules up to swapping L and R, 
this structure introduces the natural 
splitting of the whole tangent bundle $TM$, i.e. an {\em almost product 
structure} on $M$. The almost product
structure also restricts to the complex tangent bundles $T^{CR}M$.
We shall write 
\begin{gather*}
TM=T^RM\oplus T^LM,\quad \tcr M=T^{\operatorname{CR},R}M\oplus
T^{\operatorname{CR},L}M\\
\operatorname{Gr}TM=
(T^LM/T^{\operatorname{CR},L}M\oplus 
T^{\operatorname{CR},L}M)\oplus (T^RM/T^{\operatorname{CR},R}M\oplus
T^{\operatorname{CR},R}M)
\end{gather*}
but we keep in mind that the splitting is available only locally, in general.
In particular, the two components of $TM$ are orthogonal with respect to the 
algebraic bracket $\{\ ,\ \}_{\operatorname{Lie}}$. 

Next, we observe that the canonical almost complex structure $J$ 
defined on $T^{CR}M$ is induced by the $(\fg, P)$-structure. 
Indeed, we define
$$
J\in (T^{CR}M)^*\otimes T^{CR}M,\quad
J(Tp.\omi(X)(u))=Tp.(\omi(iX)(u))
$$
and this formula does not depend on the choice of $u$ and $X$ because the
adjoint action of $P$ on $\fg_{-1}\subset \fg/\ff p$ is complex linear.

At the same time, there is the obvious integrable complex 
structure coming from the definition $T^{CR}M=TM\cap i{TM}\subset T{\Bbb
C}^4$ on the embedded CR-manifolds. 
The fundamental osculation (\ref{osculation}) then implies that
these two almost complex structures on $T^{CR}M$ coincide. 

\subsection{}\label{abstract-hyperbolic}
The {\em abstract hyperbolic CR-manifolds} of dimension six and
CR-codimension two are defined by the specification of a regular $(\ff
g,P)$-structure on $M$ in the sense of Definition \ref{gp-structure}. In
particular, they come equipped by the CR-subbundle
$T^{\operatorname{CR}}M\subset TM$ of real codimension two with an almost
complex structure, and the compatible almost product structure on $TM$. 
The general theory then applies as well and so the normal Cartan connections
are given uniquely on all such manifolds. We shall see, however, that the
embedded ones have very specific features. The automatic integrability of
the almost complex structure $J$ on $T^{CR}M$ is an example. We can meet
these more general structures on some 6-dimensional real submanifolds in
8-dimensional almost complex manifolds. 

\smallskip
Our goal is to understand fully the local geometrical properties. For
that reason we shall first discuss all possible algebraic brackets on $TM$
which arise from the Lie bracket of vector fields and we shall link them to
certain components of the curvature of the canonical Cartan connection $\om$
on $M$. In fact we shall work on the abstract level, forgetting more or less
about the embedding of the manifold $M$ into ${\Bbb C}^4$. Though some of the
obstructions will vanish automatically for the embedded hyperbolic
CR-manifolds.

For example, the algebraic Lie bracket of two vector fields
$\xi, \et$ in $T^{\operatorname{CR},L}M$ has no contribution in
$T^RM/T^{CR,R}M$ and so the projection of the Lie bracket $[\xi,\et]$ to 
$T^RM=TM/T^LM$ has values in
$T^{\operatorname{CR},R}M$. Analogously we
can deal with left and right components exchanged and so there are two
obvious algebraic brackets 
\begin{gather}\label{alg-0L}
\{\ ,\ \}_L:
T^{\operatorname{CR},L}M \x T^{\operatorname{CR},L}M\to
T^{\operatorname{CR},R}M\\\label{alg-0R} \{\ ,\ \}_R:
T^{\operatorname{CR},R}M \x T^{\operatorname{CR},R}M\to
T^{\operatorname{CR},L}M 
\end{gather} 
which have to vanish automatically for all embedded hyperbolic CR-manifolds
in view of Lemma \ref{algebraic-brackets}.
We shall see in a moment that these brackets vanish even for the abstract
structures.
 
Our general strategy will be to link algebraic brackets to certain
components of the curvature $\ka$ of the Cartan
connection $\om$.  According to Theorem \ref{A-curvature}, 
we have to start by the description of the real
cohomologies 
$$ 
H^2_*(\fg_-^L\oplus \fg_-^R,\fsu(2,1)^L\oplus \fsu(2,1)^R)
.$$ 

\begin{lemma} All irreducible components of these real cohomologies are the
one-dimensional $\fg_0$-modules which are generated by the (real) bilinear
cochains listed in Table \ref{table-su}.
\end{lemma}

\begin{figure}
\centerline{\offinterlineskip
\vbox{\hrule\halign{\vrule \vrule width0pt height3.5ex depth1ex\ \hfil$#$\hfil\ 
\vrule&
\ \hfil$#$\hfil\ \vrule&\ #\hfil\vrule\cr
\text{homog.}& \text{cochains}& \text{comment}\cr
\noalign{\hrule}
1
&
\fgr_{-2}\x\fgr_{-1}\to \fgl_{-2}
&\text{real linear in both arguments}\cr
1
&
\fgl_{-2}\x\fgl_{-1}\to \fgr_{-2}
&\text{real linear in both arguments}\cr
\noalign{\hrule}
1
&
\fgl_{-1}\x\fgr_{-1}\to \fgl_{-1}
&
antilinear in both arguments
\cr
1
&
\fgl_{-1}\x\fgr_{-1}\to \fgr_{-1}
&
sesquilinear
\cr
1
&
\fgr_{-1}\x\fgl_{-1}\to \fgr_{-1}
&
antilinear in both arguments
\cr
1
&
\fgr_{-1}\x\fgl_{-1}\to \fgl_{-1}
&
sesquilinear
\cr
\noalign{\hrule}
4
&
\fgl_{-2}\x\fgl_{-1}\to\fgl_1
&
real and complex linear
\cr
4
&
\fgr_{-2}\x\fgr_{-1}\to\fgr_1
&
real and complex linear
\cr
}\hrule}}
\caption{Real cohomologies of $\fg_-$ with coefficients in $\fg$}
\label{table-su}
\end{figure}

\begin{proof}
Let us consider the $\fg_0$-modules
$$
A_\ell=H^2_{\ell}(\fg_-^L\oplus \fg_-^R,\fsu(2,1)^L\oplus \fsu(2,1)^R).
$$
By the general theory we know that the complexifications $(A^*_\ell)_\bbc$ 
of the dual $\fg_0$-modules $A^*_\ell$ are the complex cohomologies
$H^2_{-\ell}(\ff p_+,\fsl(3,\bbc)\oplus \fsl(3,\bbc))$
listed in the table of all complex cohomologies, see Table 
\ref{table-components} in Appendix \ref{B}. Further, let us notice
that the two components in $\fg_{-1}$ have a canonical complex structure.
Now, we have just to
keep in mind, that a complexification of a real linear mapping $\phi:V\to
W$, defined on a complex vector space $V$, 
splits into two components according to the splitting
of the complexification $V_\bbc=V\oplus\bar V$. If the target of such a
mapping is complex as well, then the mapping $\phi$ itself splits into the
complex linear and complex antilinear parts. Thus the complex cohomologies
on the list of Table \ref{table-components}, and the other half of them, 
must come exactly from the components listed in Table \ref{table-su}.
\end{proof}

Now we are ready to find the geometric meaning of the individual torsion
components. First, we shall focus on the  obstructions
against the integrability of the natural almost product structure on $M$. 

Thus we are interested in brackets $\operatorname{Gr}T^LM\x
\operatorname{Gr}T^LM\to \operatorname{Gr}T^RM$ and those with
the left and right components exchanged. The
restriction of $\{\ ,\ \}_{\text{Lie}}$ vanishes clearly. Hence, 
apart from the algebraic brackets
(\ref{alg-0L}), (\ref{alg-0R}), there is another candidate
\begin{gather}\label{alg-1L}
\{\ ,\ \}_L: T^LM/T^{\operatorname{CR},L}M\x T^{\operatorname{CR},L}M\to
T^RM/T^{\operatorname{CR},R}M
\\\label{alg-1R}
\{\ ,\ \}_R: T^RM/T^{\operatorname{CR},R}M\x T^{\operatorname{CR},R}M\to
T^LM/T^{\operatorname{CR},L}M
\mbox{\rlap{.}}\end{gather}
Indeed, choosing any representative of the argument from the quotient space, 
the ordinary Lie bracket projected to the desired component yields
our algebraic bracket. In contrast to the Levi form, 
these two algebraic brackets are not 
coming from the quadric by the osculation. 

\begin{lemma}\label{torsion-product0}
The brackets (\ref{alg-0L}), (\ref{alg-0R}) vanish identically. The brackets
(\ref{alg-1L}), (\ref{alg-1R}) are given by the formulae
\begin{gather}\label{alg-1Lf}
\{\pi_L(\xi),\et\}_L = -\pi_R(
Tp.\om^{-1}(\ka^{(1)}(u)(X,Y))(u))
\\\label{alg-1Rf}
\{\pi_R(\xi),\et\}_R = -\pi_L
(Tp.\om^{-1}(\ka^{(1)}(u)(X,Y))(u))
\end{gather}
where $u\in \cg$, $\pi_L$ and $\pi_R$ are the obvious quotient projections
in the left and right components of the graded tangent space,  and 
$X\in \fg_{-2}^L, Y\in\fg_{-1}^L$, or $X\in \fg_{-2}^R, Y\in\fg_{-1}^R$,
respectively, and
$$
\xi=Tp.\om^{-1}(X)(u),\quad \et=Tp.\om^{-1}(Y)(u)
.$$ 
There are no more non-trivial algebraic brackets  $\operatorname{Gr}T^LM\x \operatorname{Gr}T^LM\to \operatorname{Gr}T^RM$ and 
$\operatorname{Gr}T^RM\x \operatorname{Gr}T^RM\to \operatorname{Gr}T^LM$.
\end{lemma}

\begin{proof} We shall discuss only brackets
$\operatorname{Gr}T^LM\x \operatorname{Gr}T^LM\to \operatorname{Gr}T^RM$.
The other ones are treated analogously.

The first part is quite easy. Let us consider $\xi_x,\et_x\in
T^{\operatorname{CR},L}_xM$. Further, choose $u\in \cg$, $x=p(u)$, and
$X,Y\in\fg_{-1}^L$ such that $\xi_x=Tp.\om^{-1}(X)(u)$,
$\et_x=Tp.\om^{-1}(Y)(u)$. According to the Lemma \ref{A-brackets-lemma},
there are the projectable vector fields $\tilde\xi$, $\tilde \et$ on $\cg$
such that their projections
$\xi=Tp\o\tilde\xi$, $\et=Tp.\tilde\et$ satisfy $\xi(x)=\xi_x$,
$\et(x)=\et_x$ and 
\begin{equation}\label{proof1}
[\xi,\et](x) = Tp.[\tilde \xi, \tilde \et](u)= Tp.[\om^{-1}(X),\om^{-1}(Y)](u)
.\end{equation}
Now let us recall the general formulae (\ref{curvature1}) and
(\ref{curvature2}) for the evaluations of curvatures and remember  
there are no
curvature components of non-positive homogeneities. In particular, 
$$
\om([\om^{-1}(X),\om^{-1}(Y)](u))\in \fg_{-2}^L\oplus (\fg_{-1}^L\oplus
\fg_{-1}^R) \ \operatorname{mod}\ff p
.$$
Thus applying the projection $\pi_R$ onto the image $T^{\operatorname{CR},R}M$
of $\{\ ,\ \}_L$, we may rewrite (\ref{proof1}) as
\begin{align*}
\{\xi_x,\et_x\}_L&=\pi_R\o Tp.[\tilde\xi,\tilde\et]\\
&=\pi_R\o Tp.\om^{-1}(u)\bigl(\om([\om^{-1}(X),\om^{-1}(Y)](u))\bigr)
\\
&=-\pi_R\o Tp.\om^{-1}(u)(\ka^{(1)}(u)(X,Y))
.\end{align*} 
In particular, the bracket must vanish because there is no cohomology
represented by cochains $\fg_{-1}^L\x\fg_{-1}^L\to \fg_{-1}^R$, see Table
\ref{table-su}, and so this component of the curvature vanishes by 
Theorem \ref{A-curvature}.

We shall proceed analogously in the case of the bracket (\ref{alg-1L}). 
Let us fix again a frame $u\in \cg$, $x=p(u)$, choose 
the element in $T^LM/T^{\operatorname{CR},L}M$ represented 
by $Tp.\om^{-1}(X)(u)$ with $X\in\fgl_{-2}$, and choose another vector
$\et_x\in T^{\operatorname{CR},L}_xM$, $\et_x=Tp.\om^{-1}(Y)(u)$, with $Y\in
\fgl_{-1}$.  Next, we consider the projectable vector fields $\tilde \xi$ on
$\cg$ such that $\om^{-1}(X)=\tilde\xi$ on the image of $\si_u$ and
similarly for $\et$. Then the value of
$Tp\o\tilde\xi= \xi$ at $x$ represents the right argument in $T^LM/T^{CR,L}M$
and we obtain
\begin{equation}\label{rel-fields}
[\xi,\et](x)=Tp.[\tilde\xi,\tilde\et](u)
=Tp.[\om^{-1}(X),\om^{-1}(Y)](u)
\end{equation}
(see again Lemma \ref{A-brackets-lemma}).
Since $X\in \fg_{-2}^L$, $Y\in \fg_{-1}^L$, our table of cohomologies implies 
$$
\ka^{(1)}(u)(X,Y)=-\om^R_{-2}([\om^{-1}(X),
\om^{-1}(Y)](u))\in \fg_{-2}^R 
$$
where $\om^R_{-2}$ is the component of $\om$ valued in $\fgr_{-2}$.
In particular we obtain the required equality (\ref{alg-1Lf}). 

There are still two more possibilities for algebraic brackets
$\operatorname{Gr}T^LM\x \operatorname{Gr}T^LM\to \operatorname{Gr}T^RM$.
The first one, 
$$
T^LM/T^{\operatorname{CR},L}M\x T^LM/T^{\operatorname{CR},L}M\to 
\operatorname{Gr}T^RM
$$ 
is obviously zero since the arguments are from an one-dimensional space. 
The remaining brackets  
\begin{gather}\label{alg-2L}
\{\ ,\ \}_L: T^LM/T^{\operatorname{CR},L}M \x T^{\operatorname{CR},L}M\to 
T^{\operatorname{CR},R}M
\\\label{alg-2R}
\{\ ,\ \}_R: T^LM/T^{\operatorname{CR},R}M \x T^{\operatorname{CR},R}M\to 
T^{\operatorname{CR},L}M
\end{gather}
can be well defined and are algebraic if and only if the brackets (\ref{alg-1L})
and (\ref{alg-1R}) vanish, respectively. If so, then their values are again
defined by considering the representatives of the elements in the 
quotient spaces in the domain. By the vanishing assumption, 
their projection to the quotient on the right hand side is zero,
thus they lie in the desired targets. 

So let us assume that the bracket (\ref{alg-1L})
vanishes. Then tracing the above computation of the latter bracket 
step by step, with the target replaced by 
$T^{\operatorname{CR},R}M$, we end up with the formula
$$
\{\pi_L(\xi(x)),\et(x)\}_L=
-\pi_R\o Tp(\om^{-1}(u)(\ka^{(2)}(u)(X,Y)))
.$$
Thus the vanishing of our bracket is equivalent to the vanishing of the
corresponding component $\ka^{(2)}:\fgl_{-2}\x\fgl_{-1}\to \fgr_{-1}$.
Consider now the homogeneous component of degree two of the 
Bianchi identity, see (\ref{A-Bianchi}) in Section \ref{A}. Its right hand
side includes terms of two kinds: 
\begin{equation}\label{use-Bianchi}
\ka^{(1)}(\ka^{(1)}(X,Y),Z)\qquad 
\cl_{\om^{-1}(Z)}\ka^{(2+|Z|)}(X,Y).
\end{equation} 
The differential
$\partial\ka^{(2)}$ on the left hand side is homogeneous of degree two
again.
Since our component of $\ka^{(2)}$
is not in the list of the available cohomologies and $\ka^{(2)}$ is
co-closed, this component must be in the image of $\partial^*$. Further, 
let us notice that $\partial$ acts
injectively on the image of $\partial^*$ (cf. the Hodge-structure mentioned
in \ref{Hodge-structure}). 
Thus the image of $\ka^{(2)}$ under $\partial$
vanishes if and only if this component vanishes too. Now, we are interested
only in the component
$\fgl_{-2}{}^*\otimes \fgl_{-1}{}^*\otimes \fgr_{-1}$ and so its image under
$\partial$ will sit in the subspace (cf. (\ref{partial}))
$$
(\fgr_{-1}{}^*\otimes\fgl_{-2}{}^*\otimes\fgl_{-1}{}^*\otimes\fgr_{-2})
\oplus
(\fgl_{-1}{}^*\otimes\fgl_{-1}{}^*\otimes\fgl_{-1}{}^*\otimes\fgr_{-1})
.$$
Our knowledge of all possibly non-zero components of first degree in $\ka$
(remember we assume that the bracket (\ref{alg-1Lf}) vanishes)
and a straightforward inspection of the few possibilities of the placement of
the arguments in the two terms in the Bianchi identity shows that there is
no way to get anything non-zero. 

Thus the vanishing of the last possible algebraic bracket has been proved.
\end{proof}

\begin{theorem}\label{torsion-product1}  Let $M$ be an abstract hyperbolic
6-dimensional CR-manifold of CR-codimension two.
The left distribution $T^LM$ is involutive if and only if the bracket
(\ref{alg-1L}) vanishes, the right distribution is involutive if and only if 
the bracket (\ref{alg-1R}) vanishes. 

The almost product structure on $M$ is integrable if and only if both 
these brackets vanish. 
\end{theorem}
\begin{proof} All projections of the Lie brackets $T^LM\x T^LM\to T^RM$ are
linear over functions and thus algebraic. Therefore, Lemma
\ref{torsion-product0} implies immediately the first claim. Similarly for
the other distribution $T^RM$ and the last claim follows
by the standard foliation theory.
\end{proof}

\begin{theorem}\label{torsion-J1}
Let $M$ be a 6-dimensional abstract hyperbolic CR-ma\-ni\-fold of 
CR-codimension two.
The canonical almost complex structure $J$ on $T^{CR}M$ is integrable if and only
if the part $\ka^{(1)}_{aa}\in C^{\infty}(\cg,\ff g_{-1}^*\wedge\ff
g_{-1}^*\otimes \ff g_{-1})$ of $\ka^{(1)}$ which is antilinear in both 
arguments vanishes.
In particular, this part of the torsion vanishes on the embedded
6-dimensional hyperbolic CR-manifolds in ${\Bbb C}^4$.
\end{theorem}
\begin{proof} By the defining properties of the regular
$(\fg,P)$-structures, the complexified CR-tangent subbundle 
$T^{CR}_{\Bbb C}M\subset T_{\Bbb C}M$ must be involutive. 
Thus the obstruction against the integrability of $J$ 
is the Nijenhuis tensor $N\in
\La^2(T^{CR}M)^*\otimes T^{CR}M$. Consequently, the theorem will be proved
once we verify the following claim: 
{\it The Nijenhuis tensor $N$, expressed by its frame form 
$\nu\in C^{\infty}(\cg, \ff g_{-1}^*\wedge\ff g_{-1}^*\otimes\ff g_{-1})$,
equals to $4\ka^{(1)}_{aa}$.}

In order to prove this, let us choose vector fields $\xi$, $\et$ in
$T^{CR}M$, a frame $u\in\cg$, $p(u)=x\in M$, and $X,Y\in\ff g_{-1}$ such that
$\xi(x)=Tp.\omi(X)(u)$, $\et(x)=Tp.\omi(Y)(u)$. We have
$$
N(\xi(x),\et(x))=[\xi,\et]-[J\xi,J\et]+J([J\xi,\et]+[\xi,J\et])
$$
and $N(Tp.\omi(X)(u),Tp.\omi(Y)(u))=Tp.\omi(\nu(u)(X,Y))$.

As before, there are projectable vector fields $\tilde \xi$, $\tilde \et$
over $\xi$ and $\et$, such that $[\omi(X),\omi(Y)](u)=[\tilde
\xi,\tilde\et](u)$ and similarly for $J\xi(x)=Tp.\omi(iX)(u)$ and 
$J\et(x)=Tp.\omi(iY)(u)$. Then we can compute
\begin{align*}
N(\xi(x)&,\et(x))=Tp.\bigl([\omi(X),\omi(Y)] - [\omi(iX),\omi(iY)] +
\\
&\omi(i\om([\omi(iX),\omi(Y)] + [\omi(X),\omi(iY)])(u))\bigr)(u)
\\
=\ &Tp.\bigl(\omi\bigl([X,Y]-[iX,iY]+i[iX,Y]+i[X,iY]+
\\
&\ka^{(1)}(X,Y)-\ka^{(1)}(iX,iY) + i \ka^{(1)}(iX,Y) + i \ka^{(1)}(X,iY)
\bigr)(u)\bigr)(u)
\\
=\ &Tp.\omi(4\ka^{(1)}_{aa}(u)(X,Y))(u)
\end{align*}
\end{proof} 

\subsection{The complexified Cartan connection}\label{complex-omega}
\rm The proof of the preceding theorem could be also done by the methods of
\ref{torsion-product0}, with the help of complexification. Indeed, the
complexification of the canonical form $\om$ is $\om_{\Bbb C}:T_{\Bbb
C}\cg\to\ff g_{\Bbb C}$ which is a complex linear automorphism on each
complex tangent space. The Lie bracket of real vector fields extends to 
the complex ones and again each choice of $u\in\cg$, $X,Y\in(\ff
g_{-})_{\Bbb C}$
allows to choose projectable complex vector fields $\tilde\xi$, $\tilde\et$
such that $[\tilde\xi,\tilde\et](u)=[\omi_{\Bbb C}(X),\omi_{\Bbb C}(Y)](u)$.
Furthermore, the expansion of $\om_{\Bbb C}([\omi_{\Bbb C}(X),\omi_{\Bbb
C}(Y)])$ into the real and imaginary parts shows that the latter expression
yields exactly the complexification $\ka_{\Bbb C}$ of the curvature.
Thus we may proceed exactly as in \ref{torsion-product1} in order to link the
component of $\ka^{(1)}_{\Bbb C}$ acting on two holomorphic vectors in the
complexification of $\ff g_{-1}$ and valued in the antiholomorphic ones, with
the obstruction against the integrability of the holomorphic tangent
subbundle in $(T^{CR}M)_{\Bbb C}$. Of course, the same applies if we swap
the holomorphic and antiholomorphic vector fields.

\subsection{The rest of the torsion}\label{torsion-rest}
Similarly, the remaining two components of the curvature obtain a nice
geometric interpretation in form of an algebraic bracket which is defined as
follows. Take a holomorphic vector field $\xi\in (T^{CR,L}M)_{\Bbb C}$,
an antiholomorphic $\et\in (T^{CR,R}M)_{\Bbb C}$ and project their Lie
bracket to the holomorphic component in $(T^{CR,R}M)_{\Bbb C}$. 
Clearly, this is
an algebraic bracket and it vanishes if and only if the corresponding
curvature component vanishes. Similarly to the involutivity of the
holomorphic and antiholomorphic bundles, this obstruction has an tensorial
interpretation $S_R\in (T^{CR,L}M)^*\otimes(T^{CR,R}M)^*\otimes T^{CR,R}M$, 
\begin{equation}\label{N2-tensor-R}
S_R(\xi,\et)=\pi_R([\xi,\et]+[J\xi,J\et]-J[J\xi,\et]+J[\xi,J\et])
.
\end{equation}
Swapping the left and right tangent bundle components, we obtain 
\begin{equation}\label{N2-tensor-L}
S_L\in 
(T^{CR,R}M)^*\otimes(T^{CR,L}M)^*\otimes T^{CR,L}M.
\end{equation}

\begin{theorem}\label{hyperbolic-torsion-free}
Let $M\subset {\Bbb C}^4$ be a Levi non-degenerate 6-dimensional 
CR-manifold of CR-codimension 2 and let $x\in M$ be a hyperbolic point. 
Then $M$ is the product of two Levi non-degenerate 3-dimensional 
CR-structures $M_1,M_2\subset {\Bbb C}^2$, locally around $x$, 
if and only if the algebraic Lie
brackets 
\begin{gather*}
\{\ ,\ \}_L: T^LM/T^{\operatorname{CR},L}M\x 
T^{\operatorname{CR},L}M\to T^RM/T^{\operatorname{CR},R}M
\\ 
\{\ ,\ \}_R: T^RM/T^{\operatorname{CR},R}M\x T^{\operatorname{CR},R}M\to
T^LM/T^{\operatorname{CR},L}M
\end{gather*} 
vanish on a neighborhood of $x$.

The abstract 6-dimensional hyperbolic
CR-manifolds of CR-co\-di\-men\-sion two are locally products of (abstract)
3-dimensional CR-manifolds of CR-co\-di\-men\-sion 1 if and only if the above
algebraic brackets, 
as well as the Nijenhuis tensor $N_J\in \La^2(T^{CR}M)^*\otimes T^{CR}M$ and 
tensors  $S_R$, $S_L$ from (\ref{N2-tensor-R}), 
(\ref{N2-tensor-L})
vanish. 
\end{theorem}
\begin{proof} All considerations are local and so we may suppose that the
whole $M$ is hyperbolic. If $M$ is a product of two 3-dimensional
CR-manifolds, then we can also consider the product $\cg\to M_1\x M_2$ 
of the corresponding
canonical Cartan bundles $\cg_1\to M_1$, $\cg_2\to M_2$ equipped with the
product $\om=\om_1\oplus\om_2$
of the corresponding  normal Cartan connections. These bundles and connections were
constructed already by Cartan in \cite{Car} and their construction is also
covered by Theorem \ref{A-main}. By definition, the new form
$\om\in\Omega^1(\cg,\ff g)$ 
has all properties of normal Cartan connections and its
curvature $\ka$ is the sum of the two curvatures $\ka_1$ and $\ka_2$ of
$\om_1$ and $\om_2$, respectively. In particular, there is no torsion
because the connections $\om_1$ and $\om_2$ are torsion free. 
Thus the four tensorial obstructions
on $M\subset {\Bbb C}^4$  have to vanish as well.   

Now, let $M$ be an (abstract) hyperbolic 6-dimensional
CR-manifold and assume that all six tensorial obstructions from our theorem
vanish globally. According to previous results, all homogeneous components 
$\ka^{(1)}$ of the torsion of the normal Cartan connection vanish globally.
Thus, according to Theorem \ref{A-curvature} and the table of the relevant
cohomologies, all homogeneous components $\ka^{(i)}$, $i\le 3$, vanish too. 
In particular, there is no torsion part in $\ka$. Let us consider next the
part $\ka^L$ of the whole curvature which is represented by cochains of the
form $\ff g_-^L\x \ff g_-^L\to \ff g^R$ and analogously $\ka^R$ with left
and right components swapped. We shall use the induction on the 
homogeneity degrees to show, that all these components vanish. Thus, assume
we have done this for homogeneity less than $j$ and consider the components
in homogeneity $j$. Since there are no cohomologies of the types in question,
the corresponding parts of $\ka^L$ and $\ka^R$ are in the image of
$\partial^*$ and so the differential $\partial$ acts on them injectively.
Thus we can apply the Bianchi identity in order to see that there is
no component which could contribute, cf. the end of the proof of Theorem 
\ref{torsion-product1}.  Consequently, both $\ka^L$ and $\ka^R$ vanish.

The splitting $\ff g=\ff g^L\oplus \ff g^R$
induces two complementary $P$-invariant distributions on 
$\cg$, $T\cg=T^L\cg\oplus T^R\cg$. These
distributions are involutive if and only if the obvious algebraic bracket
$T^L\cg\x T^L\cg \to T^R\cg\simeq T\cg/T^L\cg$ vanishes and similarly with
$L$ and $R$ swapped. Since the brackets are algebraic, we may use the
parallel fields $\omi(X)$, $\omi(Y)$ with properly chosen $X,Y$ for their
evaluation. The projection may be realized by means of the component of $\om$
valued in the left or right part of $\ff g$. But this is controlled by the
curvatures $\ka^L$, $\ka^R$ and so the brackets vanish, as proclaimed.

Now, we know that the Cartan bundle $\cg$ locally splits into a product of two
manifolds but we need much more. 
We wish to prove that there is a neighborhood of $x$ over which 
the whole Cartan bundle $(\cg,\om)$ is isomorphic to a product of
$(\cg_L,\alpha_L)$ and $(\cg_R,\alpha_R)$ for some suitable Cartan connections
$\alpha_L$, $\alpha_R$. In fact, if we construct these data only locally around
a chosen frame $u\in\cg$, then the right invariance will ensure what we
need. The normal coordinates determined by the choice of $u$ will be again 
our basic tool. 

So let $\ph_u:\fgl_-\oplus\fgr_-\to \cg$ be the mapping defined only 
locally around the origin by the horizontal flows and let $\si_u$ be the
corresponding section of $\cg\to M$. By abuse of notation, we shall not
mention the definition domains of these and other locally defined mappings.
Let us write $P_L$ and $P_R$ for the parabolic subgroups in the individual
components of $G$ and define the trivial principal bundles 
$\cg_L=\fgl_-\x P_L$, $\cg_R=\fgr_-\x P_R$. Further, consider the principal
fibre bundle morphism $\Psi:\cg_L\x\cg_R\to \cg$ (notice $P_L$ and $P_R$
commute and the whole mapping is defined on fibres over a neighborhood of
the origin in $\fg_-$ only)
$$
\Psi: ((X,p),(Y,q))\mapsto \ph_u(X,Y)pq
.$$
Furthermore, the restrictions of $\Psi$ yield principal fibre bundle 
morphisms
\begin{gather*}
\Psi_L:\cg_L\to \cg,\quad (X,p)\mapsto \ph_u(X,0)p
\\
\Psi_R:\cg_R\to \cg,\quad (Y,q)\mapsto \ph_u(0,Y)q
\end{gather*}
and consider the one forms $\al_L=\Psi_L^*\om_L$, $\al_R=\Psi_R^*\om_R$ 
where $\omega_L$ and $\omega_R$ are the left and right components of
$\omega$. It remains to prove 
that $(\cg_L,\al_L)\x(\cg_R,\al_R)$ is a bundle with Cartan
connection (defined locally over a neighborhood of the origin) and $\Psi^*\om =
\al_L\oplus \al_R$ wherever defined.

First notice that, due to our choices and the involutivity of the left and
right parts of $T\cg$, the forms $\al_L$ and $\al_R$ are pullbacks of the
whole $\om$ (viewed then as forms with values in $\fg$, but
without any contribution to one half of the image). Thus the properties of
the Cartan connections are simply transfered by $\Psi_L$ and $\Psi_R$.
Furthermore, since the curvature of $\om$ does not mix left and right sides
either, the structure equations for $\al_L$ and $\al_R$ are obtained as
pullbacks of the structure equation of $\om$. In particular, the curvatures
are again $\partial^*$ closed. Thus $(\cg_L,\al_L)$ and $(\cg_R,\al_R)$ are
3-dimensional CR-manifolds of CR-codimension one 
(locally around the origin of the base manifolds). 
Finally, we observe that $\Psi^*\om$ will (locally) 
coincide with the product of the
newly constructed Cartan connections if and only if they will evaluate equally
on vectors tangent to a fixed section of $\cg_L\x \cg_R$. Thus consider the
section $(X,Y)\mapsto ((X,e),(Y,e))$, evaluate
$(\al_L\oplus\al_R)$ at the vector $(W,0)+(0,Z)\in
T_{((X,e),(Y,e))}(\cg_L\x\cg_R)$, and compare this with $\Psi^*\om$. In
fact, we may even deal with the left and right components of the tangent
space separately. 

Each such vector 
$\xi=\sfrac\partial{\partial t}_{|0}((X+tW,e),(Y,e))$ is mapped by $\Psi_L$
to $T(\Psi_L)(\xi)=
\sfrac\partial{\partial t}_{|0}\operatorname{Fl}^{X+tW}_1(u)$
and so we can easily compare the values $\al_L(\xi)$ and $\Psi^*\om(\xi)$: 
\begin{align*}
(\al_L)(\xi) &= \om_L(\sfrac\partial{\partial t}_{|0}
\operatorname{Fl}^{\omi(X+tW)}_1(u))\\
\Psi^*\om(\xi) &= \om_L(\sfrac\partial{\partial t}_{|0}
\operatorname{Fl}^{\omi(X+tW+Y)}_1(u))
.\end{align*}
Next, we observe that $\omi(Y)$ commutes with $\omi(X+tW)$ since there is no
cohomology mixing the arguments from the left and right components of
$\fg_-$. Thus we may
rewrite the last expression as 
\begin{align*}
\Psi^*\om(\xi)&=\om_L\bigl(T(\operatorname{Fl}^{\omi(Y)}_1)(
\sfrac\partial{\partial t}_{|0}
\operatorname{Fl}^{\omi(X+tW)}_1(u))\bigr)\\
&= (\operatorname{Fl}^{\omi(Y)}_1)^*\om_L
(\sfrac\partial{\partial t}_{|0}
\operatorname{Fl}^{\omi(X+tW)}_1(u))
.\end{align*}
Thus, in order to see that the two values coincide, it suffices to show that
$(\operatorname{Fl}^{\omi(Y)}_1)^*\om_L=\om_L$ for all $Y\in\fgr_-$.

We know this for the flow in the time zero,
$\operatorname{Fl}^{\omi(Y)}_0=\operatorname{id}_{\cg}$, and so we have just
to show that $\sfrac\partial{\partial s} 
(\operatorname{Fl}^{\omi(Y)}_s)^*\om_L$ vanishes identically. Each vector
in the left component of $T_v\cg$ is of the form $\omi(V)(v)$ with
$V\in\fgl$ and we compute
\begin{align*}
(\operatorname{Fl}^{\omi(Y)}_s)^*\om_L(\sfrac\partial{\partial t}_{|0} 
\operatorname{Fl}^{\omi(V)}_t(v))&=\om_L(\sfrac\partial{\partial t}_{|0} 
\operatorname{Fl}^{\omi(Y)}_s\o\operatorname{Fl}^{\omi(V)}_t(v))
\\
&=\om_L(\omi(V)(\operatorname{Fl}^{\omi(Y)}_s(v)))=V.
\end{align*}
Since the derivative of this constant mapping vanishes, the required 
invariance of $\om_L$ has been proved. 

Similarly we deal with the other component $\om_R$.

Finally we observe that if $M$ is embedded in ${\Bbb C}^4$, then 
we may always find  embeddings $\phi_i$ of the components $M_i$ in 
neighborhoods of $x_i$ into $\C^2$ such that
$$ \phi=\phi_1\oplus\phi_2:\,  M \to \C^4=\C^2\oplus\C^2$$
is an embedding of $M$ at $x=(x_1,x_2)$. In fact, consider the initial
embedding $\psi:\, M\to\C^4$. Then the restriction of $\psi$ to 
$M_1\times\{x_2\}$ is an embedding of $M_1$ into $\C^4$ that respects the
CR-structure of $M_1$. There is a holomorphic projection 
$\chi_1:\,\C^4 \to \C^2$ that is diffeomorphic from $\psi(M_1\times\{x_2\})$ 
onto its image. Denote the resulting mapping by $\phi_1$ and the analogous
mapping for the second component by  $\phi_2$. Then $\phi_1 \oplus\phi_2$
is the desired embedding. By passing to normal forms 
(see Appendix \ref{C}) one can even prove
that the embeddings $\phi$ and $\psi$ are equivalent, i.e., $\phi=\Phi\o\psi$
with some locally defined biholomorphic map $\Phi: \C^4\to \C^4$.     
\end{proof}

\subsection{Chains}
The last topic we want to discuss are the analogies to the chains on
CR-manifolds of CR-codimension one. We have introduced the general concepts
of chains and 1-chains in \ref{chains-general} for all parabolic geometries. 
These two notions coincide for the CR-manifolds of CR-codimension
one and they also coincide with the chains defined in \cite{CheM}.

Let us recall that the amount of different 1-chains up to parametrizations 
passing in fixed direction through a given point $x\in M$, 
as well as the set of all chains through $x$ is 
visible from the homogeneous case (see \ref{chains-structure}). 

The detailed discussion on the quadric $Q$ is reviewed in Appendix \ref{C}
with the following result: There is a one-parametric family of distinguished
parametrizations on each (non-parametrized) 1-chain, and in each direction
which does not belong to the subspace $T^{CR}_xM$ and does not belong to
$T^{L}_xM$ neither to $T^{R}_xM$, there is a 1-parametric class of 1-chains
up to their parameterizations, cf. (\ref{hyp-chains}).  If the direction
does belong to the left or right tangent space then there is a unique
1-chain in that direction. The chains through a given point $x\in M$ are
available only in 2-dimensional directions of the form $\{u,X\wedge Y\}\in
T_xM\wedge T_xM$ with $u\in \cg$ in the fibre over $x$ and $X$, $Y\in
\fg_{-2}$.

A general 2-dimensional surface is said to have
the {\em chain property at its point $y$} if there is a chain providing a
parametrization of this surface around $y$. 

The vector fields $\omi(X)$, $X\in \fg_{-2}$ 
span a two-dimensional distribution in $T\cg$ which
we call the {\em chain distribution} of the CR-structure.

In general, the two-dimensional (non-parameterized) chains 
$\be^{u_t}$ rotate around one fixed 1-chain 
$\al^{u_0,X}(t)$ if we move the ruling frame
$u_t=\operatorname{Fl}^{\omi(X)}_t(u)$ along the horizontal flow. This is not
possible, however, if the whole torsion of our
CR-structure is zero, because then the whole chain distribution is 
integrable. This is in accordance with the previous
theorem claiming that the whole bundle $\cg$ is the product of two canonical
Cartan bundles and the Cartan connection is a product, too. Thus, in the
torsion-free case, our chains
$\be^u$ are obtained as products of the chains in the
three-dimensional CR-manifolds. In particular we have proved the following
theorem. 

\begin{theorem}
Let $M\subset {\Bbb C}^4$ be an embedded 6-dimensional hyperbolic 
CR-manifold and assume that  the algebraic brackets
$\{\ ,\ \}_L$ and $\{\ ,\ \}_R$ from Theorem \ref{hyperbolic-torsion-free}
vanish identically. Then each chain $\be^u:U\subset
\fg_{-2}\to M$ has the chain property at each of its points. 

The same conclusion holds for abstract 6-dimensional hyperbolic CR-manifolds
of CR-codimension 2 without torsion.\end{theorem}   

\section{The elliptic structures}\label{elliptic}

\subsection{Almost complex and almost product structures}
Let us recall that on embedded elliptic 6-dimensional CR-manifolds of
CR-codi\-men\-sion two, the fundamental osculation
(\ref{osculation}) provides the $(\fg,P)$ structure on $M$ with
$\fg=\ff s\ff l(3,{\Bbb C})$, and its standard complex Borel subalgebra $\ff p$
(both viewed as real Lie algebras). The proper choices for the groups $G$, 
$G_0$, $P$ are discussed in Appendix \ref{C}.

There are striking general similarities between the hyperbolic and elliptic
geometries. Indeed,
the decomposition of the subspace $\fg_{-1}\subset \fg_-$
$$
\fg_{-1}=\fgl_{-1}\oplus\fgr_{-1}
$$ 
induces an {\em almost product structure} on the complex tangent bundle
$T^{CR}M$. We shall write again $T^{CR,L}M$ and $T^{CR,R}M$ for the
individual components. Furthermore, the
complex structure of the whole real Lie algebra $\fsl(3,{\Bbb C})$
induces the almost complex structure $J$ on $TM$, given by the formula 
$J(Tp.\omi(X)(u))=Tp.\omi(iX)(u)$. Clearly this formula is independent of
the choice of $X$ and $u$ which give the same vector $\omi(X)(u)\in T_xM$
because the adjoint action of $P$ on $\fg$ is complex linear.

As we have seen in the hyperbolic case, the knowledge of the real 
second cohomologies of the algebras in question is most essential. Also now,
we shall mostly deal with the abstract $(\fg, P)$-structures defined on
6-dimensional manifolds but we shall point out the specific properties of
the embedded ones. In particular, all obstructions coming from cohomologies
with cochains of the form $\fg_{-1}\x\fg_{-1}\to \fg_{-1}$ will disappear
automatically according to Lemma \ref{algebraic-brackets}. 

Roughly speaking, the role of the integrability of the almost complex
structures on the complex subbundles on hyperbolic manifolds is played by
the integrability of the almost product structure on $T^{CR}M$ in the
elliptic case. In particular the almost product structure will always be
integrable on the embedded elliptic CR-manifolds.  Further, the
integrability of the almost product structure of the hyperbolic manifolds
corresponds to the integrability of the almost complex structure $J$ on the
elliptic ones. In particular, the almost complex structure $J$ is intrinsic
to the manifold $M$ and it cannot be induced by the ambient complex
structure in ${\Bbb C}^4$.

\begin{lemma}
All irreducible components in $H^2_*(\fg_-,\fsl(3,{\Bbb C}))$ are the one
dimensional $G_0$-modules which are generated by the cochains 
listed in Table \ref{table-sl}.
\end{lemma}

\begin{figure}
\centerline{\offinterlineskip
\vbox{\hrule\halign{\vrule \vrule width0pt height3.5ex depth1ex\ \hfil$#$\hfil\ 
\vrule&
\ \hfil$#$\hfil\ \vrule&\ #\hfil\vrule\cr
\text{homog.}& \text{cochains}& \text{comment}\cr
\noalign{\hrule}
1
&
\fg_{-2}\x\fgl_{-1}\to \fg_{-2}
&\text{antilinear in both arguments}\cr
1
&
\fg_{-2}\x\fgr_{-1}\to \fg_{-2}
&\text{antilinear in both arguments}\cr
\noalign{\hrule}
1
&
\fgl_{-1}\x\fgl_{-1}\to \fgr_{-1}
&
sesquilinear
\cr
1
&
\fgr_{-1}\x\fgr_{-1}\to \fgl_{-1}
&
sesquilinear
\cr
1
&
\fgr_{-1}\x\fgl_{-1}\to \fgl_{-1}
&
sesquilinear
\cr
1
&
\fgl_{-1}\x\fgr_{-1}\to \fgr_{-1}
&
sesquilinear
\cr
\noalign{\hrule}
4
&
\fg_{-2}\x\fgl_{-1}\to\fgl_1
&
complex linear in both arguments
\cr
4
&
\fg_{-2}\x\fgr_{-1}\to\fgr_1
&
complex linear in both arguments
\cr
}\hrule}}
\caption{Real second cohomologies of $\fg_-$ with coefficients in
$\fg=\fsl(3,{\Bbb C})$}
\label{table-sl}
\end{figure}

\begin{proof}
Exactly as in the hyperbolic case, the complexification of the cohomologies
we want to describe is fully described by Table \ref{table-components}
in Appendix \ref{B}. Because of the complex structure on $\fg_{-}$, each of
the real components will produce two copies in the complexification. In
order to recognize them, we have to notice that complexifications of 
complex linear maps will not swap the two copies in the complexified Lie
algebra, while the antilinear ones will swap them. This simple observation
leads immediately to our Table \ref{table-sl}.
\end{proof}

\begin{theorem}\label{torsion-J2}
The almost complex structure $J$ on an abstract elliptic CR-manifold of
CR-codimension two is integrable if and only if the antilinear part
$\ka^{(1)}_{aa}$ of the curvature $\ka$ of the canonical normal Cartan
connection vanishes. This in turn happens if and only if 
the algebraic Lie brackets 
\begin{gather}\label{alg-31}
T^{(1,0)}M/T^{CR}_{\Bbb C}M\x (T^{CR,L}M)^{(1,0)} \to 
T^{(0,1)}M/T^{CR}_{\Bbb C}M
\\
\label{alg-32}
T^{(1,0)}M/T^{CR}_{\Bbb C}M\x (T^{CR,R}M)^{(1,0)} \to 
T^{(0,1)}M/T^{CR}_{\Bbb C}M
\end{gather}
on the complexified graded tangent bundles vanish
identically.
\end{theorem}
\begin{proof}
Essentially, all technique we need has been developed already. In particular,
we may repeat the computation of the Nijenhuis tensor from the proof of
Theorem \ref{torsion-J1}. Since $\fg_-$ is complex, we can do
that with any $X,Y\in\fg_-$. The result tells us that the Nijenhuis
tensor $N$, evaluated on $Tp.\omi(X)(u)$, $Tp.\omi(Y)(u)$, 
is equal to
$$
Tp.\omi\bigl((4\ka^{(1)}_{aa}+4\ka^{(2)}_{aa}+4\ka^{(3)}_{aa})(X,Y)(u)\bigr)(u)
.$$ 
Now, under the additional
condition that the higher homogeneities cannot contribute whenever
$\ka^{(1)}_{aa}(X,Y)$ vanishes, the Nijenhuis tensor vanishes if and only if
$\ka^{(1)}_{aa}$ vanishes. According to the table of
cohomologies, the latter expression must be given by the algebraic brackets
(\ref{alg-31}), (\ref{alg-32}).

Thus we have to show, that if $\ka^{(1)}_{aa}$ vanished, then no other
antilinear component valued in $\fg_-$ could occur in 
$\ka^{(2)}$, and if so, than even not in $\ka^{(3)}$.
Let us assume the two brackets (\ref{alg-31}), (\ref{alg-32}) vanish.
Then there is the algebraic bracket
$$
T^{(1,0)}M/T^{CR}_{\Bbb C}M\x (T^{CR,L}M)^{(1,0)} \to 
(T^{CR}M)^{(0,1)}
$$
which can be evaluated by means of the complexified curvature component of
homogeneity two. Clearly this must come from an antilinear component and the
vanishing of this algebraic bracket 
is equivalent to the vanishing of the antilinear parts $\ka^{(2)}_{aa}:
\fg_{-2}\x \fgl_{-1}\to \fg_{-1}$. Using the Bianchi identity exactly as in
the end of the proof of Lemma \ref{torsion-product0} we verify that
there is no curvature like this. 

Similarly we could proceed with the remaining algebraic brackets on the
holomorphic tangent bundle with values in the antiholomorhic tangent bundle.
However, the only component of homogeneity three is
$\ka^{(3)}_{aa}:\fg_{-2}\x\fg_{-2}\to \fg_{-1}$ and this vanishes
automatically because it is complex antiliear and $\fg_{-2}$ is of (complex)
dimension one.
\end{proof}

\begin{theorem}\label{torsion-product2}
Let $M$ be an abstract 6-dimensional elliptic CR-mani\-fold with
CR-codimension two. The distributions $T^{CR,L}M$, $T^{CR,R}M$ 
in the complex subspace $T^{CR}M$ are
integrable if and only if the algebraic Lie brackets 
\begin{gather}
\label{alg-33}
(T^{CR,L}M)^{(1,0)}\x (T^{CR,L}M)^{(0,1)}\to (T^{CR,R}M)^{(1,0)}
\\
(T^{CR,R}M)^{(1,0)}\x (T^{CR,R}M)^{(0,1)}\to (T^{CR,L}M)^{(1,0)}
\label{alg-34}
\end{gather}
on the complexified
complex spaces $T^{CR}_{\Bbb C}M$ vanish identically. 

In particular, these almost product structures are always integrable
on the embedded elliptic CR-manifolds.
\end{theorem}
\begin{proof}
The distributions are integrable if and only if the algebraic Lie brackets
of two fields from the same component projected to the other one vanish.
This is equivalent to the corresponding condition on the complexified
bundles $T^{CR}_{\Bbb C}$. Now we may use the technique introduced in
\ref{complex-omega}. Thus all the algebraic brackets in question will be
linked to specific components of the curvature. Since they are all living on
the CR tangent spaces, they must vanish automatically on the embedded elliptic
manifolds.

On abstract manifolds, this means the brackets of
holomorphic fields projected to the other component vanish automatically and
the distribution $T^{CR,L}M$ is integrable if and
only if the algebraic bracket (\ref{alg-33}) vanishes (cf. Table
\ref{table-sl}). The other
distribution is treated similarly.
\end{proof}

\subsection{Remaining torsion components}\label{torsion-rest2}
Let us notice that also the remaining two components of the torsion part of 
the curvature $\ka$ of the canonical normal Cartan connection allow an
expression by algebraic brackets. This time we obtain
\begin{gather}
\label{alg-35}
(T^{CR,R}M)^{(1,0)}\x (T^{CR,L}M)^{(0,1)}\to (T^{CR,L}M)^{(1,0)}
\\
(T^{CR,L}M)^{(1,0)}\x (T^{CR,R}M)^{(0,1)}\to (T^{CR,R}M)^{(1,0)}
\label{alg-36}
\end{gather}
and they vanish again on all embedded elliptic 6-dimensional CR-mani\-folds in
${\Bbb C}^4$.

\begin{theorem}\label{torsion-free2}
Let $M\subset {\Bbb C}^4$ be an embedded 6-dimensional elliptic 
CR-mani\-fold of CR-codimension 2 and assume that the algebraic brackets 
(\ref{alg-31}), (\ref{alg-32})
both vanish. Then the complex structure $J$ on the entire
Cartan bundle $\cg$ is integrable, the normal Cartan connection is
holomorphic, and there are two integrable foliations of
$M$ by complex curves in ${\Bbb C}^4$ which span the complex subbundle
$T^{CR}M$. 

The same conclusion is true on the abstract 6-dimensional elliptic
CR-manifolds if and only if all algebraic brackets 
(\ref{alg-31}), (\ref{alg-32}), (\ref{alg-33}), (\ref{alg-34}), (\ref{alg-35}),
and (\ref{alg-36}) vanish identically. 
\end{theorem}
\begin{proof}
In fact, we have nearly proved all necessary facts. Again, the same
computation with the Nijenhuis tensor reveals, that the antilinear part
$\ka_{aa}$ of the entire curvature obstructs its integrability. Once we
assume that all the torsion vanishes, there are no components of the
curvature up to homogeneity four. This is not antilinear, however. A
simple check with the Bianchi identity shows that the complex linear
curvature components can never produce anything antilinear. Thus the
integrability of the complex structure follows. Since the complex structure
$J$ on $\cg$ has been defined by the absolute parallelism $\om$, clearly
$\om\in\Omega^1(\cg,\fg)$ is holomorphic.

On the abstract manifolds, the same argument applies if we assume that the
whole torsion vanishes. On the other, each of the components of the torsion
eventually produces some antilinear contribution in higher homogeneities
via the Bianchi identity.

Now, assume $J$ is integrable and the torsion vanishes. Then also all
horizontal vector fields $\omi(X)$, $\omi(Y)$ with $X$, $Y\in\fg_{-1}^L$, or
both in the other component, commute. Thus we obtain the integrable (real)
2-dimensional distributions in $T\cg$ spanned by their values. The integral
surfaces can be locally parametrized by the holomorphic (with respect to $J$)
mappings
\begin{gather*}
\gamma^R_u:\fg_{-1}^R\to T\cg,\quad X\mapsto \operatorname{Fl}^{\omi(X)}_1(u)
\\
\gamma^L_u:\fg_{-1}^L\to T\cg,\quad X\mapsto \operatorname{Fl}^{\omi(X)}_1(u)
\end{gather*}
and also their projections to $M$ will be holomorphic curves. Obviously, we
have obtained integral manifolds for the distributions $T^{CR,R}$ and
$T^{CR,L}$.
\end{proof}

\subsection{Chains}\label{chains-elliptic}
Let $M$ be a 6-dimensional elliptic CR-manifold with CR-codimension two, 
$x\in M$, $\xi\in T_xM$. As discussed in \ref{chains-general}, the 
projections of the flows of horizontal vector fields determined by elements
in $\fg_{-2}$ are 1-chains with specific properties, while
$$
\fg_{-2}\ni X\mapsto \operatorname{Fl}^{\omi(X)}_1(u)\mapsto 
p(\operatorname{Fl}^{\omi(X)}_1(u))
$$
is the chain at $x$ determined by a fixed frame $u\in\cg$ over $x$.
A {\em complex chain} is a (locally defined) curve $\be:{\Bbb C}\to M$ which
is holomorphic with 
respect to the almost complex structure $J$ and has the chain property in
all its points.

\begin{theorem}\label{chains-torsion-free}
If the brackets (\ref{alg-31}) and (\ref{alg-32}) vanish on a neighborhood
of an elliptic point $x$ of an embedded 6-dimensional CR-manifold 
$M\subset {\Bbb C}^4$
of CR-codimension two, then there are unique complex chains 
through $x$ in all complex directions which do not belong to $T^{CR}M$. 

The same conclusion is true for the abstract elliptic CR-manifolds if 
the other four obstructions against the vanishing of the 
torsion equal to zero too.
\end{theorem}
\begin{proof} Analogously to the hyperbolic structures, there is the chain
distribution in $T\cg$ spanned by the horizontal fields $\omi(X)$ with
$X\in\fg_{-2}$. Again, the straightforward inspection of the possible
curvature components reveals that there is no curvature with
both arguments in $\fg_{-2}$ if the torsion vanishes. Thus the chain
distribution is integrable. Consequently the flows of the horizontal fields
yield  holomorphic parameterizations and the theorem is proved. 
\end{proof}

\section{Final remarks and conclusions}\label{conclusions}

\subsection{Relation to other results} Mizner \cite{Miz} constructed 
CR-invari\-ant connections for weakly uniform CR-structures of codimension 2. 
In the cases considered there, the automorphisms of the quadrics are 
always linear (thus, $\ff p_+$ is absent). Similar results were obtained by 
Garrity and Mizner for
CR-structures of codimension bigger than 2 with rigid osculating quadrics.
The CR-manifolds that are considered in this paper are not covered there.

In \cite{EIS} Ezhov, Isaev and Schmalz constructed parallelisms for
hyperbolic and elliptic manifolds. These parallelisms turn out to be
Cartan connections only in very special cases. The geometric reason for
that is the presence of torsion in our Cartan connection. We were able
to characterize the torsion-free (``semi-flat'') cases as manifolds
with integrable almost product structure in the hyperbolic case and
with integrable almost complex structure in the elliptic case. Thus we
give an answer to the question about the geometric meaning of 
``semi-flatness'' for elliptic manifolds that has been  posed 
in \cite{EIS}. 

Let us also remark that the almost CR-manifolds of CR-codimension one (e.g.
certain real hypersurfaces in almost complex manifolds) have been
studied from the point of view of the general theory of parabolic geometries
in \cite{Cap}. In particular, a nice geometric specification of the
construction from \ref{A-main} is presented there. 
 
\subsection{The parabolic CR-geometry} Unfortunately, the automor\-ph\-ism 
group of the parabolic quadric (\ref{parabolic-quadric}) does not fit into 
our scheme of general
parabolic geometries at all (notice the abuse of the non-compatible use of
the word ``parabolic'' which is used in the sense of Section \ref{A} now).
This is obvious already from its dimension which is 17. The structure of its
infinitesimal automorphisms is described in detail in \cite{Sch} and it
turns out that the discrete center ${\Bbb Z}_2$ of the hyperbolic or
elliptic group blows up into the additional dimension and one copy of
$\fsu(2,1)$ sits still inside. So it plays nicely its role of an
intermediate state between the hyperbolic and elliptic points. 

In particular the methods of Section \ref{A} which are based on the
existence of the Hodge theory on the cochains in the Lie algebra cohomology
cannot work. One should believe that some specification of the very general
approach in \cite{Mor} could be applicable. We consider this as a very
interesting open problem.

\subsection{Webster--Tanaka connections} There is a very rich underlying
geometry on each manifold equipped with a Cartan connection modelled over
graded Lie algebras. In particular, we always have the principal bundles 
$\cg_0=\cg/P_+\to M$ with structure group $G_0$ and the principal bundle
$\cg\to \cg_0$ with the structure group $P_+$. The latter bundle always
admits global smooth $G_0$-equivariant sections and the set of all of them
is parameterized by one-forms on $M$. The pullback of the $(\fg_-\oplus
\fg_0)$--component of the Cartan connection $\om$ by means of any of these
sections provides an affine connection on $TM$, i.e. a soldering form on
$\cg_0$ together with a principal connection on $\cg_0$. This construction
has been described in full generality in \cite{Sl2} and it produces exactly
the Webster-Tanaka connections on the CR-manifolds with CR-codimension one.
Thus we have a similar class of linear connections on $M$ underlying our
elliptic and hyperbolic structures.

\subsection{Natural bundles and invariant operators}
Another very interesting consequence of our construction of the canonical
Cartan connections is the theory of the semi-holonomic jet modules for
general parabolic geometries, which allows to transfer the problem of
finding invariant operators which act on some natural bundles coming from
representations of $P$ into problems in finite dimensional representation
theory. The first application of this theory is worked out in
\cite{CSS1,CSS4}. 

In particular, there are the Bernstein-Gelfand-Gelfand sequences for all
irreducible $G$-modules ${\Bbb V}$ which specialize to the BGG resolution of the
constant sheaf with coefficients in $\Bbb V$ on the homogeneous space, see
\cite{CSS4}. The analogies to classical complexes on CR-manifolds with
CR-codimension one should be localized inside of these sequences.

\renewcommand{\thesection}{\Alph{section}}
\setcounter{section}{0}

\section{Cohomologies}\label{B}

The aim of this section is to provide the list of all non-zero cohomologies
in $H^2(\fg_-, \fg)$ for the complex algebras 
$$
\ff g=\fsl(3,\bbc)\oplus\fsl(3,\bbc)\quad \ff p=\{
\text{all upper triangular matrices in $\ff g$}\} 
$$
We shall refer to the two copies of $\fsl(3,\bbc)$ as the left and right
ones. The two parts of $\ff g_0$ coincide with the parts of the Cartan
subalgebra of the diagonal matrices and all
the one-dimensional root spaces are (complex) one-dimensional. We shall
denote them as indicated in the following matrices
\begin{equation}
\fg=\begin{pmatrix}
*&\fgl_{1,0}&\fgl_{2}\\
\fgl_{-1,0}&*&\fgl_{0,1}\\
\fgl_{-2}&\fgl_{0,-1}&*
\end{pmatrix}
\oplus\begin{pmatrix}
*&\fgr_{1,0}&\fgr_{2}\\
\fgr_{-1,0}&*&\fgr_{0,1}\\
\fgr_{-2}&\fgr_{0,-1}&*
\end{pmatrix}
\end{equation} 
Here the stars fill up the subalgebra $\fg_0$, $\ff p_+$ consists of the
strictly upper triangular matrices, $\fg_1=\fg_{1,0}\oplus\fg_{0,1}$ as
$\fg_0$-module, etc.

The cohomologies for modules over simple algebras are completely described
in terms of the orbits of the Weyl groups on the weights. We shall use the
notation and technique as developed in \cite{BaE}. 
First, we have to recall a few basic facts on the representations of the 
parabolic subalgebra $\ff p\subset \fg$. 

The Dynkin diagram of $\fsl(3,\bbc)$ is 
\begin{picture}(28,10)(0,0)
\put(5,0){\hbox to0pt{\hss$\bullet$\hss}}
\put(5,2.7){\line(1,0){18}}
\put(23,0){\hbox to0pt{\hss$\bullet$\hss}}
\end{picture}. The parabolic subalgebras are denoted by crossing the nodes
which correspond to the negative simple coroots which do not belong to $\ff
p$. In our case this means one \begin{picture}(28,10)(0,0)
\put(5,0){\hbox to0pt{\hss$\times$\hss}}
\put(5,2.5){\line(1,0){18}}
\put(23,0){\hbox to0pt{\hss$\times$\hss}}
\end{picture} for both left and right 
$\fsl(3,\bbc)$. The weights of irreducible
representations of $\ff p$ are then denoted by the coefficients in 
their expressions as linear combinations of fundamental weights, placed over the corresponding
nodes. The $\fg$-dominant weights have non-negative integral coefficients,
the $\ff p$-dominant weights must be non-negative over the uncrossed nodes
only. For example, the trivial representation and the first and second 
fundamental representations of $\fsl(3,\bbc)$ have the highest weights 
$\dyn00$, $\dyn10$, $\dyn01$. Each $\ff
p$-module enjoys the filtration of the $\ff p$-submodules such that the associated
graded $\ff p$ module decomposes into the direct sum of irreducible $\ff
p$-modules. For example, the filtration and decomposition of the $\ff
p$-module $\fsl(3,\bbc)$ is as follows
\begin{equation}
\begin{matrix}
&&\fg_{-1,0}&&&&\fg_{1,0}\\
\fg_{-2}&+&\oplus&+&\fg_0&+&\oplus&+&\fg_2\\
&&\fg_{0,-1}&&&&\fg_{0,1}
\end{matrix}
\end{equation}
and in the terms of the highest weights for the one-dimensional irreducible
$\ff p$-modules
\begin{equation}
\begin{matrix}
&&\dyn1{-2}&&\dyn00&&\dyn{-1}{2}\\
\dyn{-1}{-1}&+&\oplus&+&\oplus&+&\oplus&+&\dyn11\\
&&\dyn{-2}1&&\dyn00&&\dyn2{-1}
\end{matrix}
\end{equation}

The whole Weyl group $W$ of $\fsl(3,\bbc)$ is generated by the
two simple reflections $s_1$, $s_2$ with respect to the two simple roots,
acting on the weights $\fg_0^*$.
\begin{equation}\label{B-reflections}
s_1: \dyn ab\quad \mapsto\quad \dyn {\llap{$-$}a}{a+b}\qquad s_2:\dyn ab
\quad\mapsto\quad \dyn {a+b}{-\rlap{$b$}}
\end{equation}
Since our parabolic subalgebra $\ff p\subset \fg$ 
is the Borel subalgebra, the corresponding 
parabolic subgroup $W^{\ff p}$ coincides with the whole $W$. 

The differential $\partial$ respects the homogeneities of the cochains and
so the cohomologies split into homogeneous components 
$H^*_\ell(\fg_-,\fg)$ too. Moreover, we have the identification
$H^p_\ell(\fg_-,\fg)\simeq H^p_{-\ell}(\ff p_+,\fg)$ (of
real vector spaces). Thus the {Kostant's version of the Bott-Borel-Weil
theorem} is relevant for our aims as well:
 
\begin{theorem}
Let $\ff p\subset\fg$ be a parabolic subalgebra in a complex simple algebra
$\fg$. If $A$ is a
finite dimensional irreducible $\ff g$-module of highest weight $\lambda$,
then the whole cohomology $H^*(\ff p_+,A)$ is completely reducible as a $\ff
p$-module and the 
irreducible components with
highest weight $\mu$ occur if and only if
there is an element $w\in W^{\ff p}\subset W$ such that
$\mu=w.\lambda=w(\lambda+\rho)-\rho$ and in that case it occurs in degree
$|w|$ with multiplicity one.
\end{theorem}  

The degree of an element $w\in W$ is defined as the smallest possible 
number of simple reflections whose composition is $w$. 
See e.g. \cite{Kos, Vog} for the proof the Theorem.

Now, we have a simple procedure to compute the cohomologies:
First, we write down the labelled Dynkin diagram depicting the $\ff
g$-dominant highest weight $\lambda$. For example, the highest weight of the
adjoint representation is $\dyn11$. Then we add one to each coefficient and 
act by combinations of simple reflections according to (\ref{B-reflections}). 
Finally we subtract one from each coefficient. The $\ff p$-dominant results are
just the highest weights of the cohomologies.   

Unfortunately, we deal with a sum $\fg=\fgl\oplus\fgr$ 
of two simple algebras. In order to make
use of the latter theorem, we shall view the representation spaces 
of $\fg$ and $\ff p$ as the (exterior) tensor products 
$A\btimes B$ of $\fgl$-modules $A$ and $\fgr$-modules $B$. 
In particular,  we understand the adjoint representation 
on $\fg=\fgl\oplus\fgr$ as
$$
\fg =(\fgl\btimes \bbc) \oplus (\bbc\btimes\fgr)
$$
with the obvious tensorial actions of $\ff p_+=\ff p^L_+\oplus \ff p^R_+$. 

The cohomology with values
in a direct sum of modules is just the direct sum of the cohomologies with
values in the submodules. Now, the K\"unneth theorem implies for each 
tensor product  of our modules $A\btimes B$
\begin{equation}\label{B-Kunneth}
H^p(\ff p^L_+\oplus \ff p^R_+, A\btimes B) = \sum_{i+j=p}\bigl(H^i(\ff
p^L_+,A)\btimes H^j(\ff p^R_+,B)\bigr).
\end{equation}
Thus, in order to compute the second cohomologies 
$$
H^2_\ell(\ff p^L_+\oplus\ff p^R_+,\fsl(3,\bbc)\oplus\fsl(3,\bbc))
$$ 
we have to know all
cohomologies $H^i_*(\ff p_+,\fsl(3,\bbc ))$, $i=0,1,2$, and  $H^i_*(\ff
p_+,\bbc)$. The results computed by the procedure as described above are
listed in Table \ref{table-cohomologies}.

\begin{figure}
\centerline{\offinterlineskip
\vbox{\hrule\halign{\vrule \vrule width0pt height3.5ex depth1ex\ \hfil$#$\hfil\ 
\vrule&
\ \hfil$#$\hfil\ \vrule&\ \hfil$#$\hfil\ \vrule&\ \hfil $#$\hfil\ \vrule\cr
H^0(\ff p_+,\bbc)
& 
\dyn00 
&
H^0(\ff p_+,\fsl(3,\bbc))
&
\dyn11
\cr
\noalign{\hrule}
H^1(\ff p_+,\bbc)&
\begin{matrix}
\dyn{-2}1\\
\dyn1{-2}
\end{matrix}
&
H^1(\ff p_+,\fsl(3,\bbc))
&
\begin{matrix}
\dyn{-3}3\\
\dyn3{-3}
\end{matrix}
\cr
\noalign{\hrule}
H^2(\ff p_+,\bbc)&
\begin{matrix}
\dyn0{-3}\\
\dyn{-3}0
\end{matrix}
&
H^2(\ff p_+,\fsl(3,\bbc))
&
\begin{matrix}
\dyn1{-5}\\
\dyn{-5}1
\end{matrix}
\cr
}\hrule}}
\caption{}\label{table-cohomologies}
\end{figure}

The homogeneity of the components is given by the sum of the coefficients,
which is the action by the so called grading element $E\in \fsl(3,\bbc)$,
$E=\operatorname{diag}(1,0,-1)$. There is another independent element
$F\in\fsl(3,\bbc)$, $F=\operatorname{diag}(1,-2/3,1)$, which acts trivially
on $\fg_2$, by 1 on $\fg_{-1,0}$, and by $-1$ on $\fg_{0,-1}$. Thus the action
of $F$ on a weight module is given by one third of the difference of the
coefficients over the nodes in the Dynkin diagram.

Now, the rest of our computation is quite easy since all irreducible
components in the cohomologies are one-dimensional. Thus in order to
localize the representatives of cohomologies as bilinear mappings, we have
just to evaluate the actions of the left and right $\fg_0$-elements $E^L$,
$E^R$, $F^L$, $F^R$ on the weight modules in the second cohomologies and
this always describes the possible domain and target of a bilinear
representative in the space of cochains uniquely. A half of the result is
listed in Table \ref{table-components}. The other half is obtained by
mutually replacing all the left and right components.

\begin{figure}
\centerline{\offinterlineskip
\vbox{\hrule\halign{\vrule \vrule width0pt height3.5ex depth1ex\ \hfil$#$\hfil\ 
\vrule&
\ \hfil$#$\hfil\ \vrule&\ \hfil$#$\hfil\ \vrule&\ \hfil$#$\hfil\ \vrule\cr
\multispan2\unskip\vrule\hfil\text{homogeneities}\hfil\vrule&
\multispan2 \hfil\vrule width0pt height3ex depth1ex
\text{components in $H^2(\ff p^L_+\oplus \ff p^R_+,
\fgl\otimes\bbc)$}\vrule\cr
\noalign{\hrule}
\text{total}&
\begin{matrix}
\text{actions of}\\E^L,E^R,F^L,F^R
\end{matrix}
&\text{components}& \text{cochains}\cr
\noalign{\hrule}
-1
&
2,-3,0,1
&
\dyn11\btimes\dyn0{-3}
&
\fgr_2\x\fgr_{1,0}\to \fgl_2
\cr
-1
&
2,-3,0,-1
&
\dyn11\btimes\dyn{-3}0
&
\fgr_2\x\fgr_{0,1}\to \fgl_2
\cr
\noalign{\hrule}
-1
&
0,-1,-2,-1
&
\dyn{-3}3\btimes\dyn{-2}1
&
\fgl_{0,1}\x\fgr_{0,1}\to\fgl_{1,0}
\cr
-1
&
0,-1,-2,1
&
\dyn{-3}3\btimes\dyn1{-2}
&
\fgl_{0,1}\x\fgr_{1,0}\to\fgl_{1,0}
\cr
-1
&
0,-1,2,-1
&
\dyn3{-3}\btimes\dyn{-2}1
&
\fgl_{1,0}\x\fgr_{0,1}\to\fgl_{0,1}
\cr
-1
&
0,-1,2,1
&
\dyn3{-3}\btimes\dyn1{-2}
&
\fgl_{1,0}\x\fgr_{1,0}\to\fgl_{0,1}
\cr
\noalign{\hrule}
-4
&
-4,2,0,0
&
\dyn1{-5}\btimes\dyn00
&
\fgl_2\x\fgl_{1,0}\to\fgl_{-1,0}
\cr
-4
&
-4,-2,0,0
&
\dyn{-5}1\btimes\dyn00
&
\fgl_2\x\fgl_{0,1}\to\fgl_{0,-1}
\cr}\hrule}}
\caption{}\label{table-components}
\end{figure}

In fact, we are interested in the real cohomologies $H^2_\ell(\ff
g_-,\fsu(2,1)\oplus\fsu(2,1))$ and $H^2_\ell(\ff g_-,\fsl(3,\bbc))$ where $\ff
g_-$ is the negative complement to the real Borel subalgebras $\ff p$.  As
we have mentioned already, the latter cohomologies are dual to the
real cohomologies $H^2_{-\ell}(\ff p_+,\fsu(2,1)\oplus\fsu(2,1))$ and 
$H^2_{-\ell}(\ff g_-,\fsl(3,\bbc))$. Thus the complexifications of the
requested cohomologies will be dual (as real modules) to those listed in
Table \ref{table-components}. 

\section{Normal forms}\label{C}
For embedded real-analytic hyperbolic or elliptic CR-manifolds one has
constructions of normal coordinates in the ambient space in a neighborhood
of a given point. These coordinates are uniquely determined up to 
some Lie-group action of the isotropy group of the
quadric (\ref{hyperbolic-quadric}) resp. (\ref{elliptic-quadric}). The equation
of the manifold takes then a
special form called {\em normal form} that refines the osculation
(\ref{osculation}) by
the quadric. These constructions generalize Chern--Moser's normal form
for real-analytic hypersurfaces in $\bbc^n$. They were obtained by
Loboda \cite{Lob} in the hyperbolic and by Ezhov and Schmalz \cite{ES} in
the elliptic case.

Let us recall the isotropy groups of the quadrics. It is convenient
to choose coordinates that reflect the geometric structure of the
quadrics. The hyperbolic quadric is the direct product of two hyperspheres
in $\bbc^2:$  
$$v_1=|z_1|^2, \quad v_2=|z_2|^2.$$  
The geometric structure of the elliptic quadric will be revealed by passing
to coordinates
$$w_1^\sharp=w_1+iw_2,\quad w_2^\sharp=w_1-iw_2.$$ 
Then ${\mathcal V}= v_1+iv_2=\frac{w_1^\sharp-\bar{w}_2^\sharp}{2i}$ 
and the equation of the quadric takes the form 
$${\cal V}=\frac{w_1^\sharp-\bar{w}_2^\sharp}{2i}=z_1\bar{z}_2.$$ 
Thus this quadric carries a complex structure; it is a complex
hypersurface in $\bbc^4$ with coordinates 
$z_1,\bar{z}_2,w_1^\sharp,\bar{w}_2^\sharp$. Below we will use these
coordinates and omit the sharps. 

The automorphism group of any quadric contains a transitive subgroup 
called Heisenberg group. For any point $(p,q)$ at the quadric the
Heisenberg translation that takes the origin into $(p,q)$ has the form
\begin{align*}
z^*&=z+p\\
w^*&=w+q+2i\langle z,p\rangle.
\end{align*}
 
Thus, any automorphism decomposes into  a Heisenberg translation and an
isotropic automorphism. The subgroup of isotropic automorphisms will play
the role of the parabolic subgroup $P$ with Lie algebra $\ff{p}$.   

For both our quadrics the isotropic automorphisms can be described by 
the well-known Poincar\'e formula for sphere automorphisms  
\begin{align*}
Z^*&=C(Z+AW)(1-2i\bar{A}Z-(R+iA\bar{A})W)^{-1}\\
W^*&=C\bar{C}W(1-2i\bar{A}Z-(R+iA\bar{A})W)^{-1},
\end{align*} 
where $Z$ and $W$ are diagonal $2\times 2$-matrices
with entries $z_1,z_2$ and $w_1,w_2$, respectively, 
$C$ and $A$ are complex diagonal matrices, and $R$ is a real
diagonal matrix. There occurs an additional discrete automorphism that 
interchanges $z_1 \lra z_2$ and $w_1 \lra w_2$. The group $P$ decomposes
into $G_0$ and $P_+$ where $G_0$ consists of the linear isotropic
automorphisms (with $A=R=0$, including the discrete automorphism) and 
$P_+$ consists of the ``non-linear'' automorphisms with $C=1$. The Lie
algebras of the latter subgroups are $\fg_0$ and $\ff{p}_+$, respectively.

The only difference between the hyperbolic and elliptic
case is a different definition of the complex conjugation. In the hyperbolic
case the conjugation is the usual one and $R$ is real means that it has
real entries. Thus, the automorphisms also split into a direct product.
Since the automorphism group of the sphere is $\SU(2,1)/{\Bbb Z}_3$ this shows 
that the automorphism group of the hyperbolic quadric is 
$$
\bigl((\SU(2,1)/{\Bbb Z}_3)\times (\SU(2,1)/{\Bbb Z}_3)\bigr) 
\rtimes {\Bbb Z}_2.
$$

In the elliptic case the complex conjugation is the usual one combined
with interchanging $z_1 \lra z_2$ and $w_1 \lra w_2$. $R$ is real
means now that the entries are mutually complex conjugated numbers. The 
identification of the automorphism group $G$  as  
$$
\bigl(\SL(3,{\Bbb C})/{\Bbb Z}_3\bigr)\rtimes {\Bbb Z}_2$$
is less evident than in
the hyperbolic case. As shown in \cite{Sch} the Lie algebra of infinitesimal
automorphisms of the elliptic quadric is isomorphic to $\fsl(3,{\Bbb C})$.
Since $G$ acts effectively at $\fsl(3,{\Bbb C})$ via $\operatorname{Ad}$
one can consider $G$ as a subgroup of $\operatorname{Aut}\fsl(3,{\Bbb C})$.
Both groups have the same dimension and consist of two connected components.
Therefore they must coincide. It is not hard to check that
$\operatorname{Aut}\fsl(3,{\Bbb C})\cong\SL(3,{\Bbb C})/{\Bbb Z}_3\rtimes 
{\Bbb Z}_2$.

It follows from the explicit description of the action of $G$ on the
infinitesimal automorphisms that $P$ is exactly the subgroup that respects
the filtration of $\fg$ by the $\ff{p}$-submodules (cf.
\ref{graded-algebras}).

Let us remark that the hyperbolic and elliptic quadrics have compact 
completions in $\C\bbp^2 \times \C\bbp^2$ resp. $\C\bbp^2\otimes \C$. All
automorphisms extend to automorphisms of the completion and are then 
linear with respect to the corresponding homogeneous
coordinates (see \cite{ESZ}). Moreover, these completions can be considered
as $G/P$.

Now we formulate the concrete normal form conditions: In the hyperbolic case 
the normalized equation of the manifolds takes the form

\begin{equation}\label{hyperbolic-nf}
v_j=|z_j|^2 + \sum N^j_{kl}(z,\bar{z},u),
\end{equation}
where $N^j_{kl}=\overline{N^j_{lk}}$ are  polynomials of degree $k$ in $z$
and of degree $l$ in $\bar{z}$ with coefficients that are analytic functions
of $u=\Re w$. The summation runs over all integral $k,l$ with $\max\{k,l\}>1$
and $\min\{k,l\}>0$ The polynomials satisfy the conditions  
\begin{xalignat*}{2}
\frac{\partial N^1_{k1}}{\partial \bar{z}_1}&=0 &\qquad \frac{\partial N^2_{k1}}{\partial
\bar{z}_2}&=0, \quad \text{for } k\ge 2\\
\frac{\partial^2 N^1_{21}}{\partial z_1\partial z_2}&=0 &\qquad \frac{\partial^2
N^2_{21}}{\partial z_1\partial z_2}&=0\\
\frac{\partial^4 N^1_{22}}{\partial z_1\partial z_2\partial \bar{z}_1\partial \bar{z}_2}&=0
&\qquad \frac{\partial^4 N^2_{22}}{\partial z_1\partial z_2\partial \bar{z}_1\partial
\bar{z}_2}&=0\\
\frac{\partial^4 N^1_{22}}{(\partial z_1)^2(\partial \bar{z}_1)^2}_{|u_2=0}&=0 &\qquad
\frac{\partial^4 N^2_{22}}{(\partial z_2)^2(\partial \bar{z}_2)^2}_{|u_1=0}&=0\\
\frac{\partial^5 N^1_{32}}{(\partial z_1)^3(\partial \bar{z}_1)^2}_{|u_2=0}&=0 &\qquad
\frac{\partial^5 N^2_{32}}{(\partial z_2)^3(\partial \bar{z}_2)^2}_{|u_1=0}&=0\\
\frac{\partial^6 N^1_{33}}{(\partial z_1)^3(\partial \bar{z}_1)^3}_{|u_2=0}&=0 &\qquad
\frac{\partial^6 N^2_{33}}{(\partial z_2)^3(\partial \bar{z}_2)^3}_{|u_1=0}&=0.
\end{xalignat*}
The manifold $M$ is torsion-free if and only if 
$N^j$ depends only on $z_j,\bar{z}_j$
and $u_j$. Then the conditions above coincide with Chern Moser's conditions.

In the elliptic case the normalized equation takes the form
\begin{equation}\label{elliptic-nf}
{\cal V}=z_1\bar{z}_2 + \sum N_{kl}(z,\bar{z},{\cal U}, \bar{\mathcal U}),
\end{equation}
where ${\cal U}=\frac{w_1+\bar{w}_2}{2}$ and the $N_{kl}$ are polynomials as
above (though without additional reality condition). The summation is also
as above and the polynomials satisfy the conditions 
\begin{xalignat*}{2}
\frac{\partial N_{k1}}{\partial \bar{z}_2}&=0 &\qquad \frac{\partial N_{1k}}{\partial 
z_1}&=0, \quad \text{for } k\ge 2\\
\frac{\partial^3 N_{21}}{\partial z_1\partial \bar{z}_1\partial z_2}&=0 &\qquad \frac{\partial^3 N_{12}}{\partial \bar{z}_2\partial \bar{z}_1\partial z_2}&=0\\
\frac{\partial^4 N_{22}}{\partial z_1\partial z_2\partial \bar{z}_1\partial \bar{z}_2}&=0
&\qquad \frac{\partial^4 N_{22}}{(\partial z_1)^2
(\partial \bar{z}_2)^2}_{|\bar{\eta}=0}&=0\\
\frac{\partial^5 N_{32}}{(\partial z_1)^3(\partial \bar{z}_2)^2}_{|\bar{\eta}=0}&=0 &\qquad
\frac{\partial^5 N_{23}}{(\partial z_1)^2\partial \bar{z}_1(\partial \bar{z}_2)^2}_{|\bar{\eta}=0}&=0\\
\frac{\partial^6 N_{33}}{(\partial z_1)^3(\partial\bar{z}_2)^3}_{|\bar{\eta}=0}&=0. &&
\end{xalignat*}

Elliptic torsion-free manifolds are obtained for 
$N=N(z_1,\bar{z}_2,\mathcal U)$.
Then $M$ is a complex hypersurface in $\mathbb C^4$ with complex
coordinates $(z_1,\bar{z}_2,w_1,\bar{w}_2)$.

{From} the normal form one can see that the real 2-dimensional surface
$\{z=0,\;v=0\}$ resp. $\{z=0,\;{\cal V}=0\}$ is always contained in the manifold.
It is called standard 2-chain $\mu_0$ (with respect to the given
normalization). One can define analytic 2-chains as all possible images of the
standard 2-chain under renormalizations. The family of chains passing through
a given point does not depend on the choice of normal coordinates but it
does depend on the initial point. In difference to the situation for
hypersurfaces a 2-chain $\mu$ for the initial point $p$ need not be a
2-chain for other points $p'\in \mu$.  

It is easy to obtain the analytic 2-chains for the quadrics through the origin
as the images of the standard chain under isotropic automorphisms. Thus
these 2-chains coincide with the geometrically defined chains from 
\ref{chains-general}. One obtains that
unparameterized 2-chains are the intersections of the quadric with so-called 
matrix lines $Z=AW$, where $Z,W,A$ have the same meaning as above.
Automorphisms with $C=1,\,A=0$ preserve the standard chain and change only
the parameter. Since the renormalizations of a manifold coincide up to 
higher order terms
with automorphisms of the osculating quadric this shows that there exists
exactly one 2-chain through the origin tangent to $\{Z=AU\}$. 

The 1-chains considered in \ref{chains-general} can be easily described
for quadrics (see \cite{Sch}). The projections of the 1-parametric 
families from $\fg_{-2}$
are straight lines in $\mu_0$ through the origin. All other 1-chains are
obtained by the action of isotropic automorphisms. Since the latter preserve
2-chains it follows that 1-chains always remain in some 2-chain. The
isotropic automorphisms decompose into one automorphism that preserves
$\mu_0$ and one that maps $\mu_0$ to another chain. Therefore it suffices
to study the 1-chains that are contained in $\mu_0$. For the hyperbolic
quadric we have the following situation:
\begin{itemize}
\item There are two singular directions at $\mu_0$ such that the only
1-chains in these directions are straight lines: $\{u_1=0\}$ and $\{u_2=0\}$.
\item In all non-singular directions one has a 1-parametric family of
1-chains consisting of one straight line and hyperbolas
\begin{equation}\label{hyp-chains}
u_1=\frac{\alpha u_2}{1-\beta u_2},
\end{equation}
where $\alpha$ indicates the direction and $\beta$ is the additional
parameter.
\item {2-chains} may intersect at single points or singular 1-chains only.
\end{itemize}
In the elliptic case we have
\begin{itemize}
\item In any direction of $\mu_0$ there is a 1-parametric family of 1-chains
consisting of a straight line and circles
\begin{equation}\label{ell-chains}
\beta(u_1^2+u_2^2)+\sin\alpha u_1 -\cos\alpha u_2=0,
\end{equation}
where $\alpha$ indicates the direction and $\beta$ is the additional
parameter.
\item 2-chains intersect at single points.
\end{itemize}

\bigskip
\noindent
Mathematisches Institut der Universit\"at Bonn, Beringstra\ss e 1, D-53115
Bonn, Germany
\\[2mm]
Department of Algebra and Geometry, Masaryk University, Jan\'a\v ckovo n. 2a, 
662 95 Brno, Czech Republic

\end{document}